\documentclass[english]{article}
\usepackage[T1]{fontenc}
\usepackage[latin9]{inputenc}
\usepackage{amsmath}
\usepackage{graphicx}

\usepackage{babel}

\begin{document}

\title{A Victorian Age Proof of the Four Color Theorem}

\author{I. Cahit\\\small{\texttt{email:icahit@gmail.com}}}
\date{}
\maketitle

\begin{abstract}
In this paper we have investigated some old issues concerning four color map problem.  We have given a general method for constructing counter-examples to Kempe's proof of the four color theorem and then show that all counterexamples can be rule out by re-constructing special $2$-colored two paths decomposition in the form of a double-spiral chain of the maximal planar graph.\\ In the second part of the paper we have given an algorithmic proof of the four color theorem which is based only on the coloring faces (regions) of a cubic planar maps. Our algorithmic proof has been given in three steps. The first two steps are the maximal mono-chromatic and then maximal dichromatic coloring of the faces in such a way that the resulting uncolored (white) regions of the incomplete two-colored map induce no odd-cycles so that in the (final) third step four coloring of the map has been obtained almost trivially.
\end{abstract}

\section{Introduction}

Four color map coloring problem is to color regions of a (normal) map $M$ with at most four colors so that neighbor regions (countries) would have receive different colors. This simple problem posed and conjectured to be true for all maps by Guthrie in 1852 [1],[32]. Its correct proof was first given in 1976 and repeated several times by the same method by the help of a computer [2]-[5]. The author has given two non-computer proofs of the four color theorem based on spiral chains in planar graphs [6],[7],[8].

In this paper we will give another one based on step-wise mono-chromatic coloring, two coloring and then four coloring of any given normal map $M$, i.e., four coloring of the faces of any cubic planar graph. Therefore our proof suits with the mathematics of the Victorian age [9],[33] in which the four color problem arose. In order to make a smooth transition to the proof we will re-investigate particularly counter-examples ("bad" examples) to Kempe's proof. Michael Rosellini in his undergraduate project summaries existing proofs together with the historical initial efforts. For his study of an counter-example he has chosen the paper of  Holroyd and Miller entitled "The example that Heawood should have given" [10] which is actually same example given by Errera [31] but drawn in the plane differently [11]. A close look to that example reveals a property which leads to a general method for constructing a class of counter-examples. On the otherhand we have given a method to re-color vertices of the "bad" maximal graph around the undecided degree five vertex for which Kempe's argument may fail, so that under the resulting four coloring the graph is decomposed into edge disjoint two paths. Furthermore the shape of the paths as seen from the Figure 1 is a double-spiral chain centered at the undecided vertex. Of course any four coloring of $G$ induces edge disjoint two bipartite graphs but not necessarily connected and in the form of a double-spiral. We have also suggest surveys on the early developments of the four color problem by Saaty [12] and Mitchem [13].

The notion of equitable colorability was introduced by Meyer [17]. That is the sizes of color classes differ by at most one. Similarly equitable labeling of graphs introduced by the author in 1990 [18]. However, an
earlier work of Hajnal and Szemérdi [19] showed that a graph $G$ with degree  $\Delta(G)$
is equitably $k$-colorable if $k\geq\Delta(G)+1$. In 1973, Meyer formulated the following
conjecture:\\

\textbf{Conjecture 1} \emph{(Equitable Coloring Conjecture (ECC) [17]). For any connected graph $G$, other than a complete graph or an odd cycle, $\chi_{=}(G)\leq\Delta(G)$.}\\

The Equitable $k$-Coloring Conjecture holds for some classes of graphs, e.g.,
outerplanar graphs with $\Delta \geq 3$ [20] and planar graphs with $\Delta\geq 13$ [21]. However the four colorings given for bad-examples in Figure 1 are all equitable $4$-coloring.\\

We have the following claim:\\

\emph{Claim.} Let $G$ be a maximal planar graph. Then there exits $4$-coloring of $G$ for which at least the sizes of three color classes differ by at most one.

\begin{figure}
\centering
\includegraphics[scale=0.4]{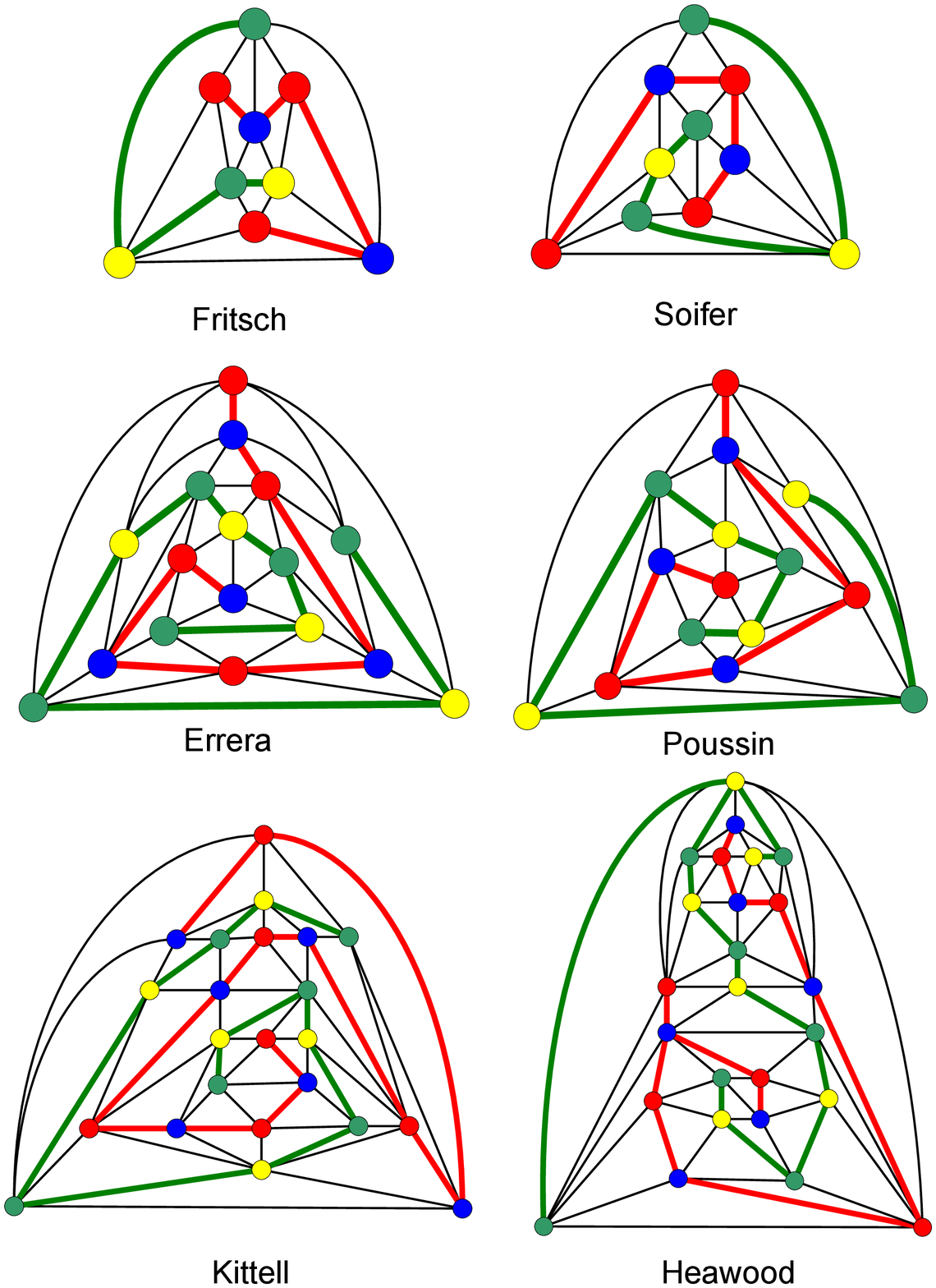}
\caption{All known counterexamples to Kempe's "proof" with double-spiral chain decompositions.}
\end{figure}
\begin{figure}
\centering
\includegraphics[scale=0.35, angle=-90]{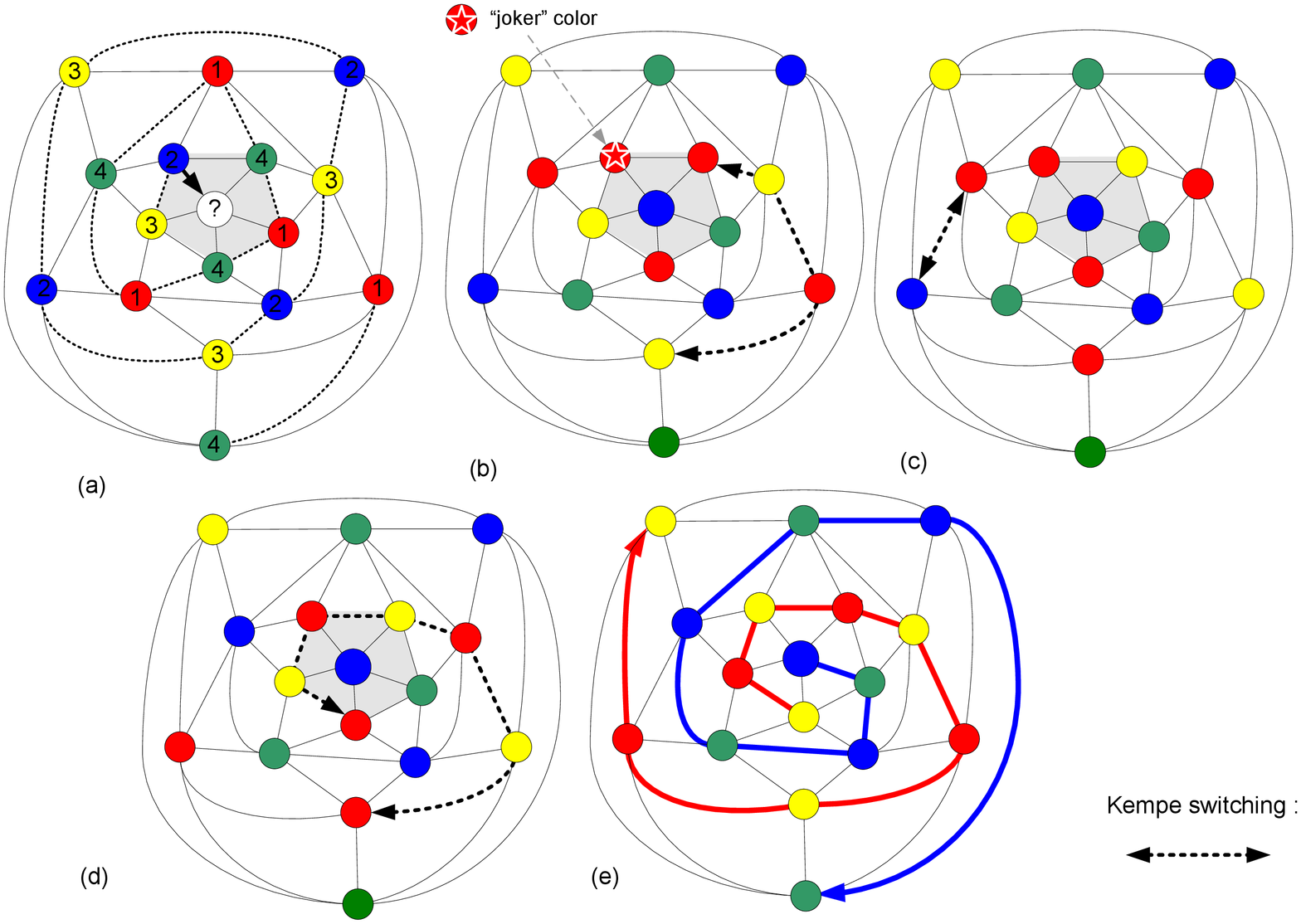}
\caption{Step-by-step resolution of an impasse in the Errera's graph.}
\end{figure}

\section{Bad Examples for Kempe's Argument}

After studying all known bad-examples to Kempe's argument one can reach to the conclusion that it is occurred only for specific planar graphs with specific incomplete four-coloring. Gethner et. al. [22],[23] have investigated Kempe's flawed proof of the Four Color Theorem from a computational and historical point of view. Kempe's "proof" gives rise to an algorithmic method of coloring planar graphs that sometimes yields a proper vertex coloring requiring four or fewer colors. They also investigate a recursive version of Kempe's method and a modified version based on the work of I. Kittell [30].

Let $G$ be an maximal planar graph with $n$ vertices.Let $T$ be the triangulation of $G$. Let $G_{1}\in \{P_{1},C_{1}\}$ and $G_{2}\in \{P_{2},C_{2}\}$ be two vertex disjoint paths or cycles such that $|G_{1}|\approx|G_{2}|$ and $|G_{1}|+|G_{2}|=n$ if under such a decomposition of $G$ every triangle $t_{i}$ has exactly one edge either from $G_{1}$ or $G_{2}$ then we say triangulation $T$ as $\beta$-triangulation. If $|C_{1}|\equiv|C_{2}|\equiv 0 (mod 2)$ then $4$-coloring of $G$ easily can be obtained. For example in Figure 1 for Fritch's graph $|P_{1}|=3,|P_{2}|=4$, for Sofier's graph $|P_{1}|=4,|P_{2}|=3$, for Errera's graph $|P_{1}|=8,|P_{2}|=7$, for Poussin's graph $|P_{1}|=7,|P_{2}|=6$, for Kittell's graph $|P_{1}|=10,|P_{2}|=11$ and finally for Heawood's graph $|P^{*}_{1}|=12,|P^{*}_{2}|=11$, where $P^{*}_{1},P^{*}_{2}$ are acyclic graphs. We choose the four colors as $\{\textbf{R}ed, \textbf{B}lue, \textbf{Y}ellow, \textbf{G}reen\}$ or for another reason $\{\textbf{B}rown, \textbf{G}reen,\textbf{d}ark \textbf{B}lue, \textbf{l}ight \textbf{B}lue\}$ or $\{1,2,3,4\}$. Moreover \emph{white} colored vertex or region in a map means awaiting color from the four-color set.

One of the important property of an "real" bad-example to Kempe's argument is that occurrence of  Kempe tangling must be independent from the order of the selection of Kempe-chains. For example Errera's bad example (first incomplete $4$-coloring of Figure 2) satisfies this condition. Now consider $C_{5,in}=\{B,G,R,G,Y\}$ that surrounds undecided white vertex. Consider also two disjoint $2$-colored cycles of length six (shown dashed lines),i.e., $C_{6,in}=\{R,G,R,G,R,G\}$ and $C_{6,out}=\{B,Y,B,Y,B,Y\}$ which forms an triangulated ring [24]. After cyclically shifting the colors in $C_{6,in}$, insert the Red "joker" color instead of Blue vertex in $C_{5,in}=\{B,G,R,G,Y\}$. Then the three Kempe chain switchings; $Ch(R,Y,R,Y),Ch(R,B)$ and $Ch(R,Y,R,Y,R,Y,R)$ (see Figure 2) resolves the impasse and a double spiral chain results [25].

\subsection{Construction of a class of bad-examples}
A triangulated ring is a $2$-connected planar graph $G_{r}$ with two faces $F_{i}$ and $F_{o}$ whose facial walks are the (induced) cycles $C_{i}$ and $C_{o}$ respectively such that: (a) $V(C_{i})\cup V(C_{o})=V(G)$ and $V(C_{i})\cap V(C_{o})=\phi$ where indices $i$ and $o$ are being used to denote the inner and outer cycles (faces)of the graph and (b) every face other than $F_{i}$ and $F_{o}$ is a triangle. We further assume that all triangles in $G_{r}$  are of type $\beta$-triangle, that is exactly one edge of the triangle belongs $C_{i}$ or $C_{o}$. Since we are interested in small size "bad-example" graphs we consider only $|C_{i}|=|C_{o}|=4,6$.
Let us give a simple lemma first.\\

\textbf{Lemma 1.} \emph{A triangulated ring $G_{r}$ with a $\beta$-triangulation and with $|C_{i}|=|C_{o}|\equiv 0 (mod 2)$ can be $4$-colored such that $C_{i}$ and $C_{o}$ colored disjoint $2$-color classes.}

\emph{Proof.} Since the inner and outer cycles are of even length; color inner cycle, say with blue and red and outer cycle with green and yellow. The $\beta$-triangulation of $G_{r}$ prevents any color conflicts in the four coloring.

Now we can construct a maximal planar graph $G$ from $G_{r}$ as follows:
(i) Place an edge $e_{i}$ inside of the inner face $F_{i}$ and place also an edge $e_{o}$ inside of the infinity (finite if the map embedded on sphere) outer-face $F_{o}$.
(ii) Make a maximal planar graph $G$ by joining the end vertices of $e_{i}$ with the vertices of $F_{i}$ and by joining the end vertices of $e_{o}$ with the vertices of $F_{o}$ such that resulting triangulation is a $\beta$-triangulation and $e_{o}$ is an outer-edge of $G$. We say inner-cycle $C_{i,in}$ is a \emph{handcuffs} for the inner-edge $e_{i}$. Similarly we say outer-cycle $C_{i,out}$ is a \emph{handcuffs} for the outer-edge $e_{o}$. The reason of this terminology will be clearer when we extract bad-examples for Kempe's argument from $G$. We will be interested in the following four coloring of $G$: Color vertices of $C_{i,in}$ and $e_{o}$ by $R$ and $B$ colors and color vertices of $C_{i,out}$ and $e_{i}$ by $Y$ and $G$ colors. This four coloring of $G$ is an proper coloring since under the cycle and edge decomposition, the triangulation is a $\beta$-triangulation. In case of cycles are of length six, let $C_{6,in}=\{u_{1},u_{2},...,u_{6}\}$, $e_{o}=\{u_{7},u_{8}\}$ and $C_{6,out}=\{v_{1},v_{2},...,v_{6}\}$, $e_{i}=\{v_{7},v_{8}\}$. Let us assume that under the $\beta$-triangulation of $G$ we also have two special Kempe-chains as follow:\\
(i) $(Y,R)$-chain $\Rightarrow$ $ch(v_{1},u_{2},v_{7},u_{6},v_{5},u_{7})$\\
(ii) $(Y,B)$-chain $\Rightarrow$ $ch(v_{1},u_{1},v_{7},u_{3},v_{5},u_{8})$\\
Now we are ready to construct the twin-bad-example graphs for Kempe's argument.\\

\emph{(a)} \textbf{Twin-graph $G_{1}$.} (\emph{Trouble in inner-face}). Delete any two edges, other than the edges of $C_{6,in}$ and $e_{i}$, of $\beta$-triangulation bounded by $C_{6,in}$ and $e_{i}$ such that the resulting new face $F_{5,in}$ contains the edge $e_{i}$ in its boundary cycle of length $5$. For example we have deleted edges $(v_{7}u_{3})$ and $(v_{8}u_{3})$ from $G$ and obtain a new cycle (face) $C_{5,in}=(v_{7},u_{2},u_{3},u_{4},v_{8})$. Now we claim that under the existing four coloring of $G$ if we place a new vertex $v_{x}$ inside of face $F_{5,in}$ and join all vertices of $C_{5,in}$ to vertex $v_{x}$ then the resulting incomplete four coloring of the modified planar graph $G_{1}$ is an bad-example to Kempe's argument. That is the four colors appear in $C_{5,in}=\{v_{7},v_{8},u_{4},u_{3},u_{2}\}$, (i.e., see Figure 3(b): $(Y,G,R,B,R)$)cannot be reduced to three colors by any Kempe-chain switching. One reason of this impasse is that $(Y,G)$ (resp. $(R,B)$) end-vertices colored edge $e_{i}$ (resp. $e_{o}$) cannot be extended due to $(R,B)$ (resp. $(Y,G)$ colored handcuffs cycle. Moreover $(G,B)$-chain $ch(v_{8},u_{5},u_{8},v_{2},u_{3})$ and $(B,Y)$-chain $ch(u_{3},v_{3},u_{8},v_{1},u_{1},v_{7})$ would prevent to reduce the number of colors to three on the vertices of $C_{5,in}$. Hence incomplete four coloring of the maximal planar graph $G_{1}$ with $17$ vertices shown in Figure 3(b) is an bad-example to Kempe's argument.\\
Note that we have the same decomposition as above if we consider;\\
$(G,R)$ cycle $C_{6,in}=\{v_{6},u_{6},v_{8},u_{4},v_{4},u_{7}\}$ and $e=\{u_{2}v_{2}\}$ and\\
$(Y,B)$ cycle $C_{6,in}=\{v_{1},u_{1},v_{7},u_{3},v_{3},u_{8}\}$ and $e=\{u_{5}v_{5}\}$.\\

\emph{(b)} \textbf{Twin-graph $G_{2}$.} (\emph{Trouble in outer-face}). The second bad-example graph $G_{2}$ can be obtained from $G$ by deleting edges $(v_{1}u_{7})$ and $(v_{1}u_{8})$. Outer-cycle of $G_{2}$ is $C_{5,out}=(u_{7},v_{6},v_{1},v_{2},u_{8})$ that has been colored by $R,G,Y,G,B$ (see Figure 3(b)). Now if we place the new vertex $v_{x}$ in the outer-face of $G_{2}$ and join to the vertices of $C_{5,out}$ then $v_{x}$ cannot be colored by the use of Kempe's argument.\\
This due to the $(Y,R)$- and $(Y,B)$-chains mentioned in (i) and (ii) before. Moreover switching of colors of the end-vertices of the edges $(v_{6}u_{1})$ or $(v_{2}u_{2})$ would not reduce the number of colors on $C_{5,out}$. Hence $G_{2}$ is an bad-example graph to Kempe's argument.

In Figure 4(a) and (b) we have shown another twins bad-example graphs $G_{1}$ and $G_{2}$ with $13$ vertices where the handcuffs cycles $C_{4,in}$ and $C_{4,out}$ are of length four. Moreover comparing the known-bad-example graphs shown in Figure 1 the graphs $G_{1}$ and $G_{2}$ are the smallest bad-examples in which  occurrence of an impasse is not depend on the order of Kempe chain switching.In Figure 4(c) we also have illustrated double-spiral chain four coloring of the bad-example of Figure 4(b). It is not difficult to show that this is possible for all bad-example graphs [24].\\

In the next section we propose a new proof for the four color theorem without using Kempe-chains based on step-by-step coloring of the faces of cubic planar maps.

\begin{figure}
\centering
\includegraphics[scale=0.4]{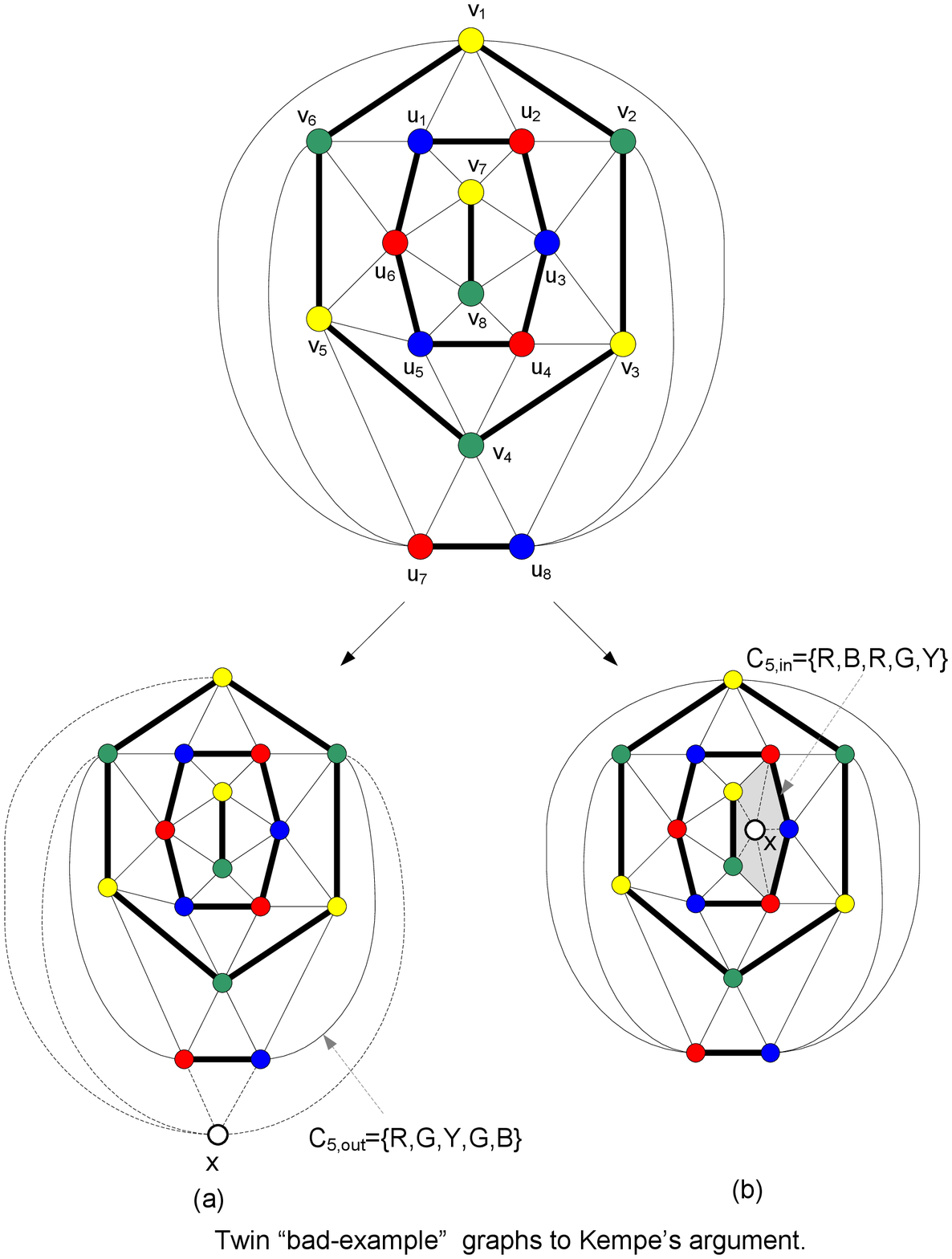}
\caption{Four coloring of an generator maximal planar graph with $\beta$-triangulation: (a) and (b) twin bad-example graphs for the Kempe's argument.}
\end{figure}
\begin{figure}
\centering
\includegraphics[scale=0.4]{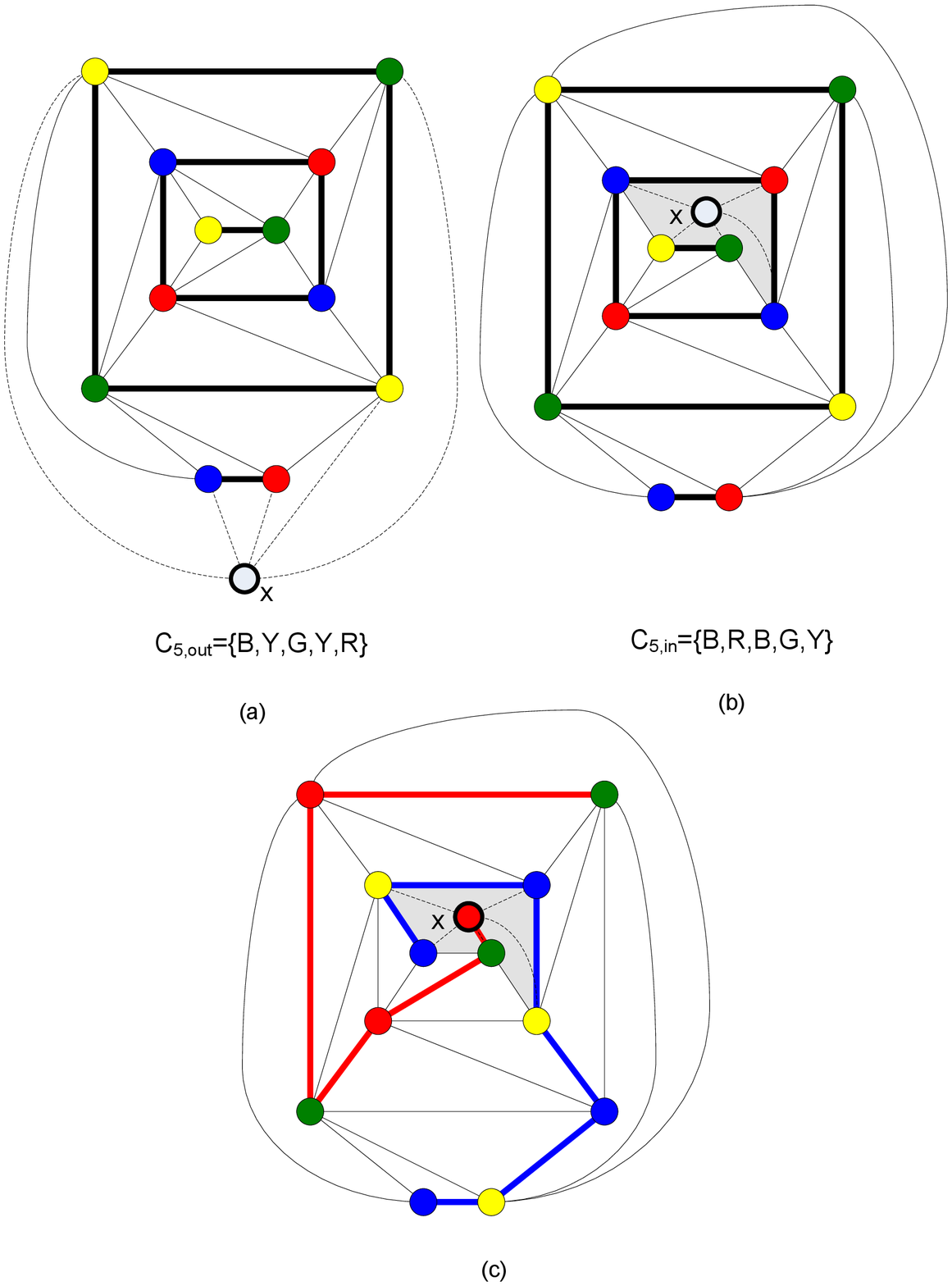}
\caption{ (a),(b)Four coloring of twin bad-example maximal planar graphs with $2C_{4}\cup \{e_{in}\}\cup \{e_{out}\}$ and (c) double-spiral chain coloring of (b)}
\end{figure}

\begin{figure}
\centering
\includegraphics[scale=0.4]{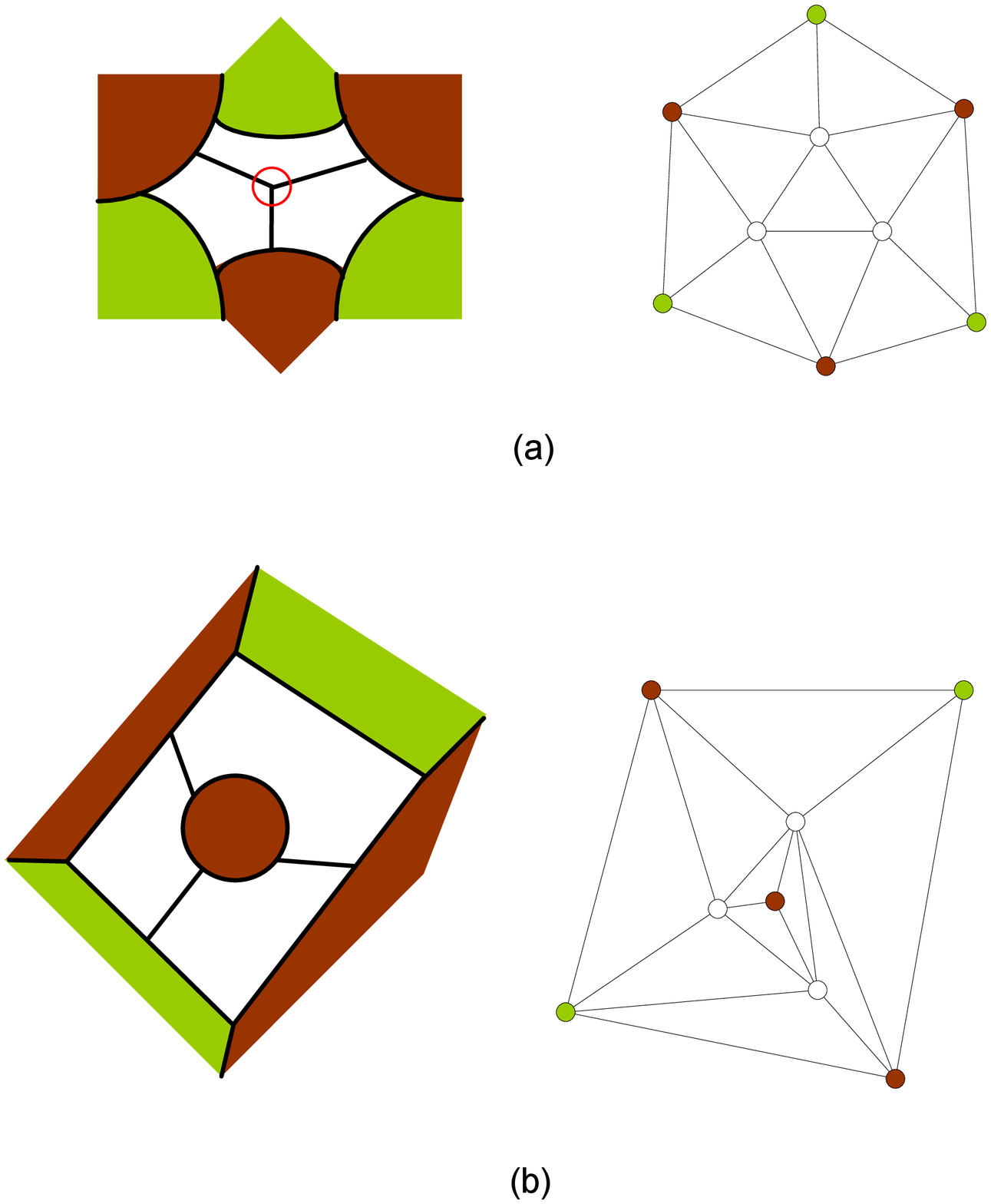}
\caption{Bad configurations in a maximal two colored map that require five colors: (a) with unwanted spot (this was the case of a bad-example to Kempe's argument; see 2-color handcuffs cycle $C_{6}$, (b) without unwanted spot but with white odd-ring.}
\end{figure}
\begin{figure}
\centering
\includegraphics[scale=0.4]{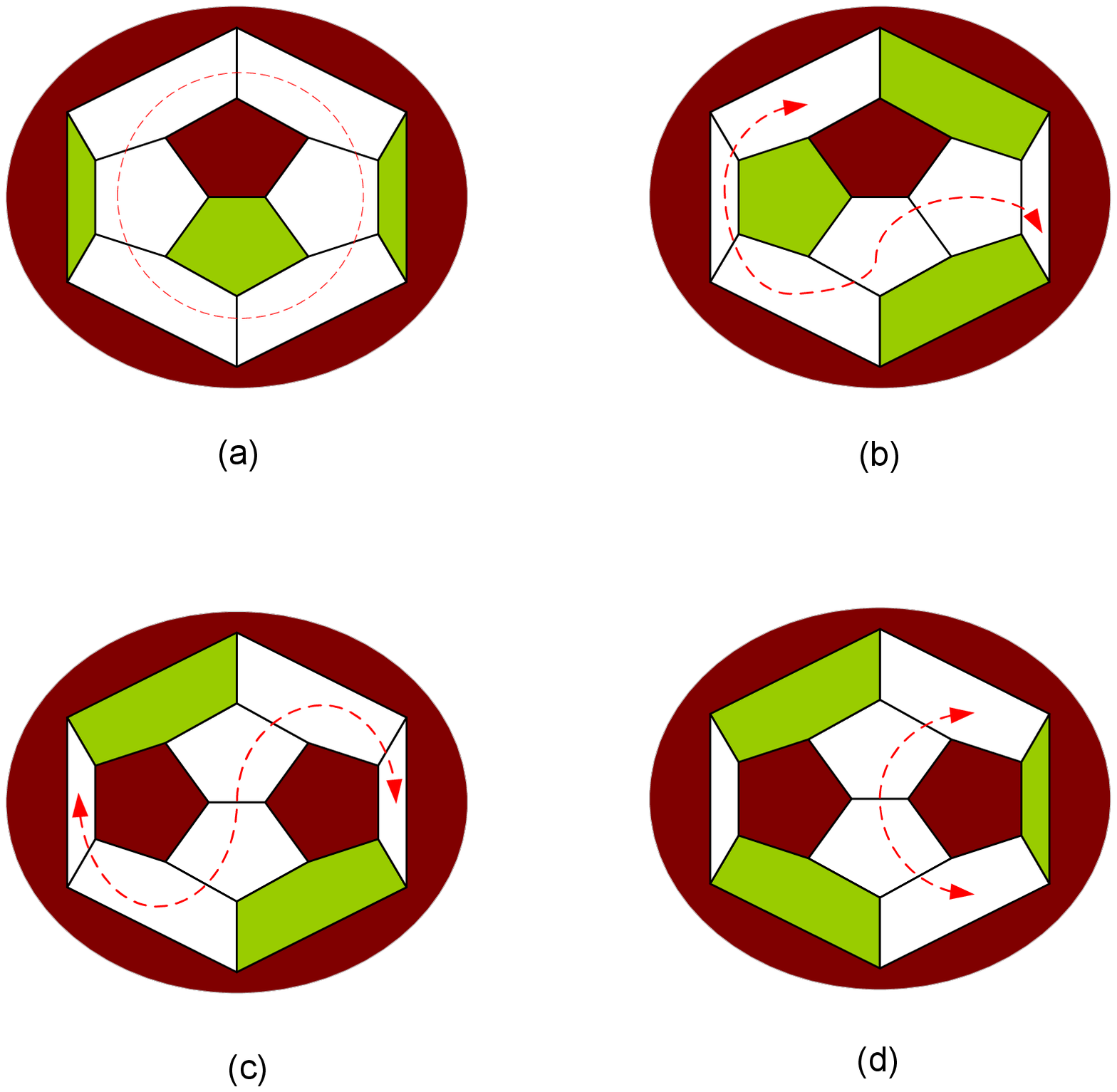}
\caption{All maximal dichromatic colorings of the Birkhoff's diamond creates no bad-configuration.}
\end{figure}

\section{A New Proof of the Four Color Map Theorem}

A more courageous title of this section would be "\emph{How to create a four colored world in three steps?}"
It is well-known and without doubt that four color theorem is true. What are the reasons for a lengthy existing proofs by the use of a computer? One answer would be going to the long way which has been forced by the false Kempe's "proof", see for example Birkhoff's reducibility of double $C_{5}$ (actually overlapped $4$ cycles of length $5$)[14]. Another answer would be over looking difficulties of the planar three colorability problem in the light of Grotzsch and Heawood's theorems [15],[16].
Starting point of the classic tedious proofs of the four color theorem is based on Kempe's failure of the reducibility of the pentagon case. That is the correct proof gives up on pentagons and turns to larger \emph{reducible} configurations for which Kempe's argument is sound [1]-[5]. The first such configuration which has ring-size $6$, was discovered by Birkhoff [14] known as Birkhoff's diamond. The Birkhoff's diamond is also inspirational starting point in our spiral ordering  and maximal mono-chromatic and maximal dichromatic coloring of the maps which leads to an non-computer proof. For example in Figure 6 we have shown all possible maximal dichromatic colorings of the Birkhoff's diamond that creates no bad configurations. That is they all be extended to proper four colorings of the map.
In this section we will be giving a new proof of the four color map theorem in which we have implicitly by pass the three-coloring problem of planar graphs within the constructive proof.

In fact our algorithmic proof implies the following theorem without relying on the four color theorem [26],[27]:\\

\textbf{Theorem 1.} \emph{Every planar graph can be decomposed into the edge disjoint union of two bipartite graphs.
}\\

Let us denote by $\emph{M}$ an normal map with $n+1$ regions, where $(n+1)$th region $r_{n+1}$ is the outer-region of $\emph{M}$. Without loss of generality we may further assume that $\emph{M}$ is digon-free (two-side region) and triangle-free (three-side region). Since if the map has a digon or triangle we shrink it to a point. Since then we can four color the resulting map and put back digons and triangles; it's surrounded by at most three colors. so there is a spare color to color the digon or triangle, as required. $\emph{M}$ can be equivalently represented by a cubic planar graph $G_{c}(\emph{M})=(V_{c},E_{c})$, where $V_{c}$ is the set of vertices associated with the crossing of pairwise three neighbor regions, and $E_{c}$ is the set of edges in the form of Jordan curve associate with the boarder of two neighbor regions between two vertices.
In order to make the map-coloring algorithm more visible and meaningful let us define the four-color set as $C=\{B,G,dB,lB\}$, where\\
- $B$ denotes \emph{brown} color and when it is assigned on to the \emph{white} background color the corresponding region becomes a "high-land".\\
- $G$ denotes \emph{green} color and when it is assigned on to the \emph{white} background color the corresponding region becomes a "low-land".\\
- $dB$ denotes \emph{dark-blue} color and when it is assigned on to the \emph{white} background color the corresponding region becomes a "deep sea".\\
- $lB$ denotes \emph{light-blue} color and when it is assigned on to the \emph{white} background color the corresponding region becomes a "shallow-sea".\\
Initially the given map colored all by background color white and at the end of the coloring algorithm (three steps) it will be colored by the colors $C$ and no white color remains on the map. Clearly we will show that this is always possible for any map $\emph{M}$.

By $M(B)$ we denote a map in which maximal number of its regions colored by $B$ (mono-chromatic coloring) where the term maximal means that any additional brown region (high-land) results color conflict and all the remaining regions are background-color white. Similarly by $M(B,G)$ we denote a map obtained from $M(B)$ in which maximal number of its white regions colored by $G$. Hence $M(B,G)$ is an maximal two-coloring of $M$.\\

\textbf{Definition 1.} \emph{In a mono-chromatic coloring of map $M(B)$ if an vertex $v$ is not incident to any  brown colored region then $v$ is called unwanted-spot or simply a spot . Furthermore if the map $M(B)$ is spot-free then the map $M(B)$ is called clean map.}\\

\textbf{Definition 2.} \emph{Spiraling of a map $M$ is a process of ordering and labeling the faces (regions), starting from the outer-region
$r_{n+1}$ and selecting always outer next region $r_{i}$ neighbor to the previous region $r_{i+1}$ in the form of a spiral.}\\

Note that depending on the adjacency of the regions of the map $M$ we may have several spirals but the ordering of the regions is
uniquely determined by the initial region and next one with the direction selected e.g., clockwise or counter clockwise.
Start with the outer face and label it $r_{n+1}$. Then, draw a curve from a
point in $r_{n+1}$, crossing an edge into an adjacent face. Label that face $r_{n}$.
The faces adjacent to $r_{n}$ (apart from $r_{n+1}$) are ordered clockwise; choose the
first (i.e. leftmost) such face, cross into it, and label it $r_{n-1}$. Proceed in this
fashion, always crossing into the leftmost available face that has not been
visited already.
At some stage one will be unable to proceed. If all faces have been visited,
then the spiral chain $S_{1} = \{r_{n+1},...,r_{1}\}$ is the spiral ordering. Otherwise,
start at the closest face to the last face of $S_{1}$ and produce a new chain $S_{2}$.
And so on till all faces are in some chain.
Similar definition has been given for maximal and cubic planar graphs in [6],[7].
For an illustration spiraling see the nested three spirals shown in blue, red and green colors in Figure 7.\\
\begin{figure}
\centering
\includegraphics[scale=0.36, angle=-90]{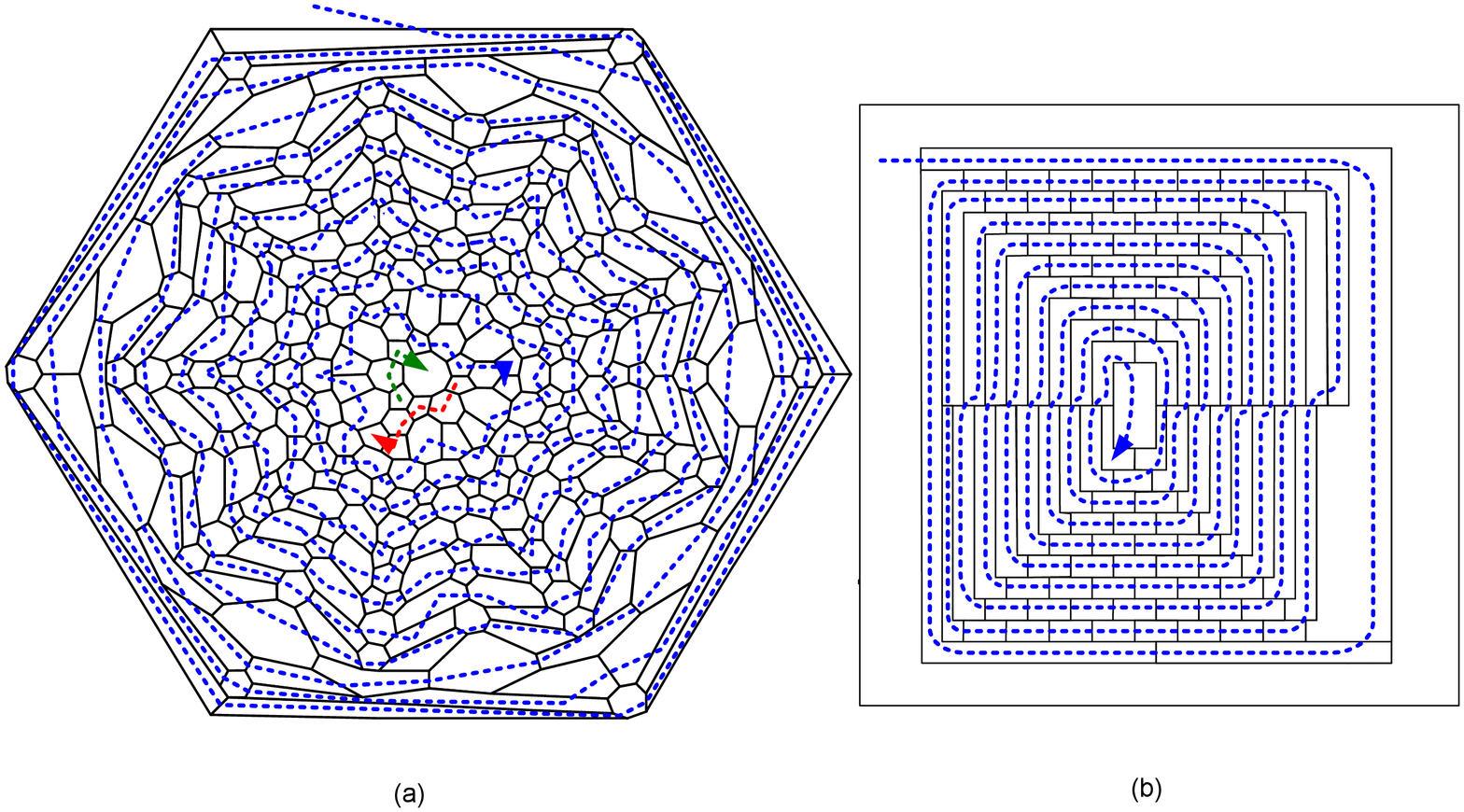}
\caption{Spiraling of the Haken and Appel's and Martin Gardner's maps.}
\end{figure}

\subsection{The map coloring algorithm.}

Main feature of the coloring algorithm is the use of each of the four colors one-by-one and preparing the conditions satisfied for the next step.\\

\emph{\textbf{Step 1.}} \emph{Maximal mono-chromatic coloring of high-lands map $M(B)$.}\\ Let $S=\{r_{n+1},r_{n},r_{n-1},...,r_{1}\}$ be the spiral ordering of the faces of map $M$. Color outer-face $r_{n+1}$ of $M$ with $B$. Along the spiral $S$ color next white region $r_{i}\in S$ with $B$ by the following rule:\\
(i) All the first neighborhood of the region $r_{i}$ remain in white (uncolored).\\
(ii) If any white region $r_{j},j>i$ is colored, that is $c(r_{j})=B$ then a color-conflicts arises.\\
(iii) At least one of the second neighborhood region $r_{i}$ with maximum number of sides would be colored by $B$. Note that all vertices of the map $M$ can be considered as spots since it is cubic planar and no region colored in brown.\\
Using (i)-(iii) and spiraling $S$ the maximal mono-chromatic set of $k$ regions can be obtained. Let us call the map $M$ after the coloring as $M(B)$.
Let us also denote the spots of $M(B)$ with a set $P=\{p_{1},p_{2},p_{3},...,p_{k}\}$ where $k<n$. That is $P$ is the set of triply neighbor white regions of the map $M(B)$ where some of the white regions may be overlapped.\\
The output of the step 1 is simply maximal disjoint of highland islands all colored in brown.\\

\emph{\textbf{Step 2.}} \emph{Maximal dichromatic coloring of high-low-lands map $M(B,G)$.}\\
We use the same spiraling $S$ of the map $M(B)$. While assigning color green $G$ to a white region consider the following two conditions:\\
(i)  Along the spiral ordering when assigning green color to white regions give priority to the white-region which has maximum number of spot vertices in $M(B)$;\\
(ii) Do not create any $(B,G)$-ring $R(B,G)$ which contains an inside odd white-ring $R(W)$ and do not leave any spot vertices.\\

We have also the following simple property of $M(B)$.\\

\textbf{Lemma 2.} \emph{The spots of the triply neighbor white regions of the map $M(B)$ cannot induces a cycle.}\\

\emph{Proof.} Let us assume that a region $r$ colored by $B$ has been surrounded by an cycle of spot vertices. Hence regions in the second neighborhood must be also all white. But (iii) we have colored at least one of the region in the second neighborhood in $B$ and that breaks the cycle of the spots into a path.\\
As it has been seen that Step 1 is rather straight forward and map $M(B)$ can easily be obtained for any $M$. Assuming the maximal mono-chromatic coloring of $M(B)$ as a base, it is not such an easy task to obtain dichromatic map $M(B,G)$. In the next step we will give the details and proofs that starting from mono-chromatic $M(B)$ it possible to two-coloring of $M(B,G)$ with a set of properties that satisfies four colorability of the whole map. That is we will show that by assigning color green (color for low-land) to the some of the white regions of $M(B)$ we obtain maximal dichromatic coloring of $M(B,G)$ without any spots, without any even $(B,G)$-ring and without odd any $W$-ring (white-rings in $M(B,G)$).\\
Let us remind the role of two-colored even cycles (handcuffs) in constructing new counter-examples to Kempe's argument in Section 2.\\
In Figure 5 we have demonstrated one of reason of an bad assignment of color green in $M(B)$. That is even-ring  $R(B,G)$ would prevent to complete coloring of white regions with four colors.\\

\textbf{Lemma 3.} \emph{Mono-chromatic (green) spiral-chain coloring of the white regions of the map $M(B)$ results in a spot-free map $M(B,G)$.}\\

\emph{Proof.} If a spot-vertex remain in $M(B,G)$ it would be one of the bad configurations illustrated in Figures 5. But this bad configuration can only occur when green color assigned without considering the maximum number of spots of the white region. However this has been protected by Step 2 (i) in the algorithm.\\

\textbf{Lemma 4.} \emph{The maximal di-chromatic map $M(B,G)$ obtained by the algorithm has no odd-white-ring of length $3$.}\\

\emph{Proof.} Lemma follows since we assumed that the cubic planar map $M$ has no triangle and $M(B,G)$ is spot-free.\\

Along the spiral ordering of $B$-coloring (brown color) of the regions of $M$ if we do \emph{not} give priority to the region with the maximum number of spots (Step 1(iii)) i.e., region with maximum sides, then there exists certain counter-examples that spiral coloring algorithm requires the fifth color. Simplest map with this property is shown in Figure 8(a) together with a dichromatic coloring of $M(B,G)$ when the face $R$ is assumed as the outerface. This is possible since the map can be redrawn so that the chosen region is the outside. The resulting spiral chains is shown in Figure 8(b). Since $M(B,G)$ of Figure 8(a) has two odd-white rings of length $5$, it cannot be extended to four coloring. It is straightforward that pentagons $P,Q$ must receive the same color in any $4$-coloring, since any attempt with $c(P)\neq c(Q)$ results in use of fifth color on the region neighbor to the outerface. This can only be resolved by the use of appropriate Kempe's switching that results in $c'(P)=c'(Q)$. On the otherhand as in Figure 8(b) if in $B$-coloring we select region $X$ with maximum number of sides instead of brown square of Figure 8(a), the coloring of $M(B,G)$ ends up without odd-white cycles.

To show that spiral ordering uses no more than four colors, consider $B$-coloring of Figure 8(c). Spiral ordering is shown in red-dashed curve. Outerface is colored by $B$ as usual. The second region colored by $B$ is the square region on the left side of the map. However Step 3 (iii) has not selected next square (colored in light blue) since the neighbor region $X$ has eight sides. The last region to be colored by $B$ is $R$ which is the first region of Figure 8(b). Now Step 2 of $G$-coloring chooses regions $Y$ and the left and right regions since they all have $8$ sides and removes the $4$ spots. Clearly $M(B,G)$ has no odd-white rings and can easily be extended to a four coloring.

\subsection{Blocking big-odd-white cycles in $M(B,G)$.}

From the above discussion it is possible to claim that map coloring algorithm would not generate $M(B,G)$ with odd-white rings (cycles). It is difficult to prove this claim, in the absence of control of detecting odd-white cycle of size greater than $3$ in the algorithm.
We will give a simple dynamic binary labeling algorithm that maintaining the parities of all white rings till the end of the map coloring algorithm. Let $M(B)$ be the maximal mono-chromatic coloring of $M$. Let $R=\{r_{n+1},r_{n},...,r_{1}\}$ be the set of regions of $M(B)$. Clearly $c(r_{n+1})=B$ but we do not certain about other regions (islands) in color $B$. Any attempt to color a white-region with $B$ results a color conflict in $M(B)$.
Let $M_{o}(B,G)$ be the set of all maximal dichromatic maps that have at least one odd-white cycles of length greater than $3$. Call these odd cycles in $M_{o}(B,G)$ as big-odd-white (simply odd-cycle) cycles.\\

Define an binary labeling $f$ of an region $r_{i}\in R$ as follows:\\

$f(r_{i})=\{ 1$   {\emph{if} $|r_{i}|\equiv 1 (mod 2)$ and $0$   {\emph{if} $|r_{i}|\equiv 0 (mod 2)\}$.}\\

Beginning of the Step 2 ($i=1$) we have the set $F_{i}$ of binary labels\\

$F_{i}= F_{i}(W) \bigcup F_{i}(B) \bigcup F_{i}(G)$\\

Initially $F_{i}(B,G)=\emptyset$, (for $i=1$) where $F_{i}(B,G)$ is the set of all binary labels corresponding to the disjoint sub-maps $M_{i}(B,G)$ computed by \\

 $f_{i}(M(B,G))=\sum_{r_{j}\in M_{i}(B,G)} f(r_{j})  (mod 2)$\\

and $F_{i}(W)=\{f(r_{j})| \emph{if}\hspace{0.2cm}   c(r_{j})=W \}$,  $F_{i}(B)= \{f(r_{j})| \emph{if} \hspace{0.2cm}  c(r_{j})=B \}$ and $F_{i}(G)= \{f(r_{j})| \emph{if} \hspace{0.2cm}  c(r_{j})=G \}$.\\

It is apparent that at each color "green" assignment of Step 2 the size of the set $F_{i}(B,G)$ increases at most by one and priority is given to the white region with maximum spots and with the even sizes. At some step $k$ we have obtained\\

$F_{i}(B,G)= \{f(M_{1}),f(M_{2}),...,f(M_{k})\}$\\

which corresponds to all disjoint sub-maps of the maximal di-chromatic coloring of $M(B,G)$. If all $f(M_{i})=0, i=1,2,...,k$ then $M(B,G)$ is white-odd ring free.\\

\textbf{Definition 3.} \emph{Let $M(B,G)$ be a maximal dichromatic map. A region $r_{j}$ colored by $B$ or $G$ is called an island if it is surrounded by the white regions.
}\\

Let $M_{s}(B,G)$ be a maximal connected two colored sub-map of $M(B,G)$. Then the white ring $R(M_{s}(B,G))$ that surrounds $M_{s}(B,G)$ is an white-odd ring if and only if the number of odd faces (including induced white faces) in $M_{s}(B,G)$ is odd. That is white ring $R(M_{s}(B,G))$ is odd iff $f_{s}M(B,G)=\equiv 1 (mod 2)$.

This property is quite useful when we are testing whether or not an odd-ring surrounds a maximal connected sub-map in $M(B,G)$.\\

We will show that by using simple transformations on odd cycles cycles it is possible to remove all odd cycles from  $M_{o}(B,G)$. This will be illustrated on the counter-example map $M(B,G)$ given in Figure 8(a) in two ways:\\

(a) \emph {Each individual odd ring blocked separately}.\\
Consider the sub-map $M_{s'}(B,G)$ containing green pentagon $Q$. Use Kempe switching for the regions of $M_{s'}(B,G)$. Now the upper and lower white regions neighbor to the outerface are all surrounded by brown-white regions. Therefore we can assign color green to these regions (denoted by $b_{1}$ and $b_{2}$) and block the odd ring that surrounds $M_{s'}(B,G)$. Now consider the sub-map $M_{s''}(B,G)$ containing the pentagon $P$. The odd ring that surrounds $M_{s''}(B,G)$ can be easily blocked by joining the two islands (brown square regions) by assigning green color to the regions denoted by $b_{3}$ and $b_{4}$ and changing the green region of $M_{s''}(B,G)$ into white color. Now as shown in Figure 9(a) the new two colored map is free of white odd rings.\\

(b)\emph {Blocking the two overlapped odd rings.}\\
Use Kempe-switching for the regions of the sub-map $M_{s'}(B,G)$ and move green regions to the white region neighbor to the pentagon $P$. Now the white regions denoted by $X$ and $Y$ which are common for both odd-rings are surrounded by brown and white regions. Therefore we can assign color green ($b_{1}$) to either region $X$ or $Y$ and block both odd rings. In Figure 9(b) odd-rings are block by re-coloring region $Y$ as $c(b_{1})=G$.\\

From the spiral coloring algorithm as well as blocking big-odd-white cycles we observe that if the number of disjoint sub-maps $M_{s}(B,G)$ is minimum then there is no odd white cycles in the map $M$.

The following lemma is useful and confirms the discussions above in general.\\

\textbf{Lemma 5.} \emph{Let $M(B,G)$ be a maximal dichromatic map which has minimum number of disjoint sub-maps $M_{s}(B,G)$. Then there is no odd-island in $M(B,G)$.}\\

\emph{Proof.} Since the regions neighbor to odd-island must be visited by the spiral ordering if the odd-island colored by $B$ (brown) in the first round of spiral coloring at least one neighbor region must be colored by $G$ (green) in the second round of spiral coloring.\\

Note that there is a map for which spiral coloring produces an even-island (see Figure 16 square island labeled $6$ in the Tutte's map). Brown squares in Figure 8(a) are not counted as island since they can be joined to the (motherland) sub-map containing pentagon $P$ without violating maximality of dichromatic coloring.\\

The following lemma is important in establishing the non-computer proof of the four color theorem. Enables to break any odd white ring in the dichromatic coloring of the map by re-arranging green colored regions  around the brown regions.\\

\textbf{Lemma 6.} \emph{Let $M(B,G)$ be a maximal di-chromatic map without spots. Let $R(W)$ be a white odd-ring of length greater than $3$ in $M(B,G)$. Then $M(B,G)$ can be re-colored to make $M'(B,G)$ white odd-ring $R(W)$ free.
}\\

\emph{Proof.} In fact we will show that for a given maximal di-chromatic coloring of $M(B,G)$ by the use of certain Kempe-chains switching it is possible to re-color any selected white region $r(w)$ (odd or even sizes) into green $G$ or brown $B$ color. Let $r(w)\in R(W)$ where $R(W)$ is denoting white odd ring of size greater than $3$. Simplest case is when region $r(w)$ is an square $S$. Let $\{r_{1},r_{2},r_{3},r_{4}\}$ be the regions neighbor to $S$. The two non-adjacent regions must be white, say $c(r_{2})=c(r_{4})=W$ and if the other two non-adjacent regions are green $c(r_{1})=c(r_{3})=G$ then color $c(r(w))=B$ or if $c(r_{1})=c(r_{3})=B$ then color $c(r(w))=G$. If $c(r_{1})=B$ and $c(r_{3})=G$ then we use Kempe switching to the $(B,G)$-Kempe chain where $r_{1}$ is the first region in the chain. So we can assign $c(r(w))=B$. Or we use Kempe switching for the $(G,B)$-Kempe chain, where $r_{3}$ is the first region in the chain. So we can assign $c(r(w))=G$. Hence the case of white square $S$ is settled.
Now consider the case when $|r(w)|\geq 5$ which is a bit different than the above. Let $r(w)$ be a pentagon $P$ with neighbor regions $\{r_{1},r_{2},r_{3},r_{4},r_{5}\}$. Clearly if $c(r_{1})=c(r_{3})=G$ and $c(r_{2})=c(r_{4})=c(r_{5})=W$ then we assign $c(P)=B$. So we assume $c(r_{4})=B$ (or $c(r_{5})=B$) but not $c(r_{2})=G$ since the white pentagon $P\in R(W)$. Now consider a $(B,W)$-Kempe chain in $M(B,G)$ starting from the region $r_{4}$ such that;\\
\emph{(i) }$(B,W)$-Kempe chain ends up with a region colored by $B$ surrounded by all $W$ and $G$ regions or ends up with a region colored $W$ surrounded by all $W$ and $G$ regions and \\
\emph{(ii)} when Kempe switching applies to $(B,W)$-Kempe chain no spots have been generated.\\
Hence after the Kempe switching of $(B,W)$-Kempe chain we have $c(r_{1})=c(r_{3})=G$ and $c(r_{2})=c(r_{4})=c(r_{5})=W$ and we assign again $c(P)=B$.
Now if $|r(w)|>5$ we apply one by one $(B,W)$-Kempe chain switching  for each brown neighbor region of $r(w)\in R(W)$ and make all neighbors of $r(w)$ all colored $G$ and $W$. Then by coloring $c(r(w))=B$ we block the white odd-ring $R(W)$. Note that no Kempe's tangling has occurred here since $R(W)$ divides the dichromatic map $M(B,G)$ into two parts and Kempe switching applied one at a time.\\

In Figure 12 although there is no white odd ring we have demonstrated by the use of the Kempe-chains how the regions labeled by $X$ and $X'$ can be recolored by brown by freeing up all neighbor regions from the color brown.\\

\textbf{Theorem 2.} \emph{The map $M(B,G)$ obtained by the Map-Coloring-Algorithm in Step 2 can be extended to a four coloring of $M$.}\\

\emph{Proof.} Proof follows from Lemmas 2,3,4 and 6.\\

\begin{figure}
\centering
\includegraphics[scale=0.3]{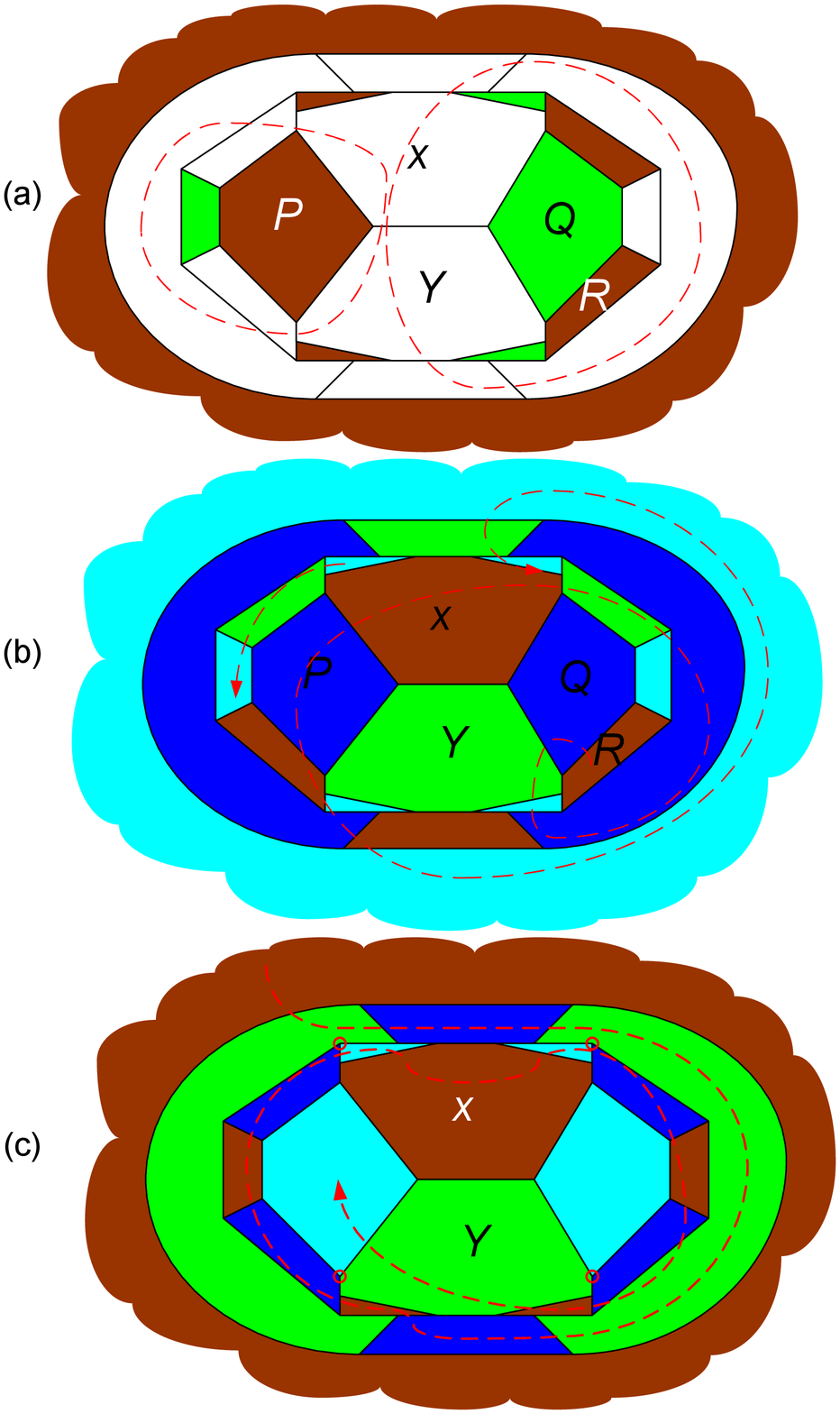}
\caption{(a) A counter-example, (b) its resolutions when $R$ is chosen as an outerface and (c) when using spiral ordering from the outerface.}
\end{figure}

\begin{figure}
\centering
\includegraphics[scale=0.35]{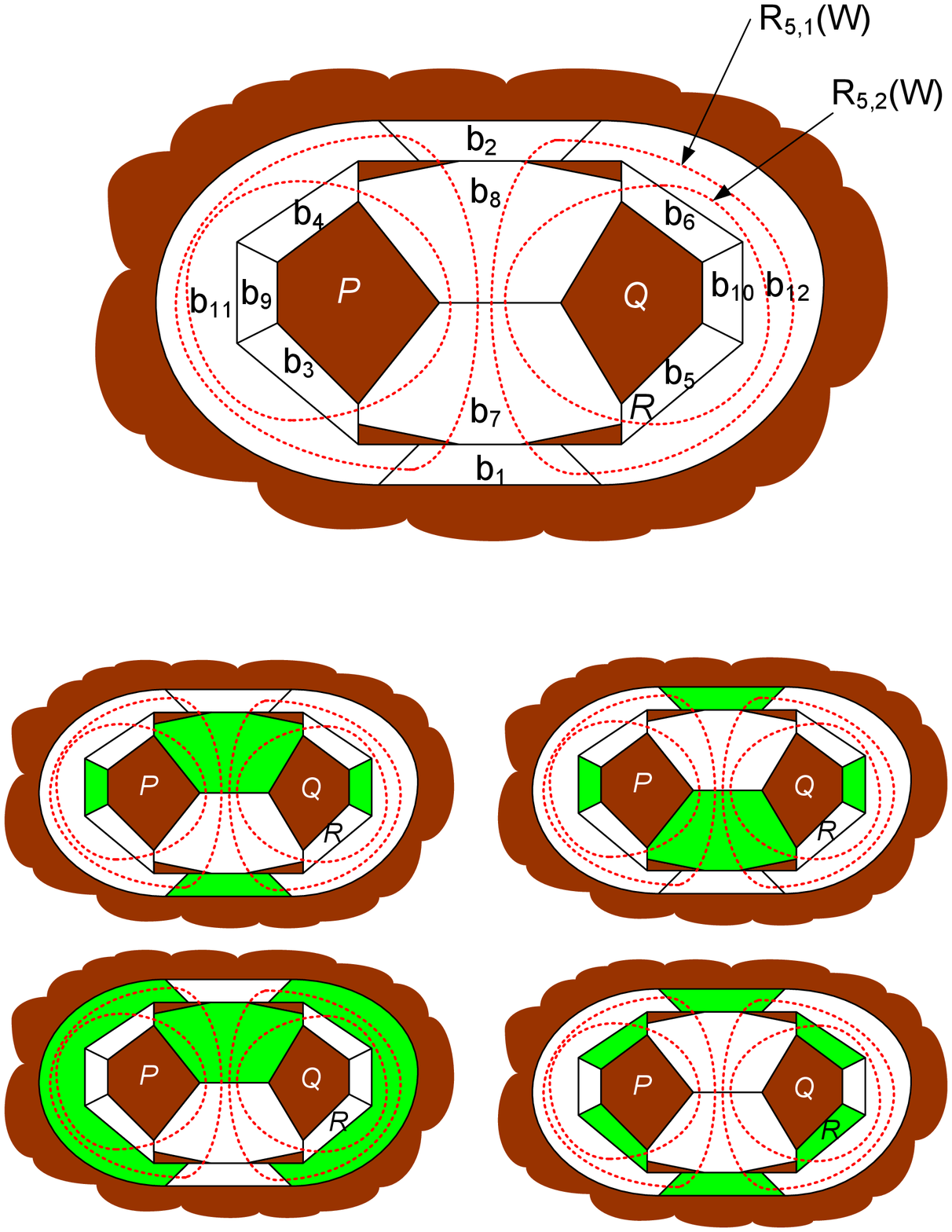}
\caption{Blocking big-odd-white rings}
\end{figure}

\emph{\textbf{Step 3.}} Four coloring of $M(B,G,lB,dB)$.\\
Since maximal dichromatic map $M(B,G)$ has only even white-rings and acyclic white regions, i.e., forest of disjoint trees and paths we can easily color them with light-blue $lB$ and dark-blue $dB$.

That is at the end of Step 3 the initial all-white normal map $M$ transformed into four colored map of $M(B,G,lB,dB)$ with the regions of high-lands, low-lands, deep-seas and shallow-seas.\\

From Theorem 2 we re-state the famous four color map theorem.\\

\textbf{Theorem 3.} \emph{All cubic planar maps are $4$-colorable.}

\vspace{1cm}
\subsection{Some illustrations of the map coloring algorithm}
The map coloring algorithm has been illustrated by the two well-known maps. Figure captions give the details.
A normal map $M$ has a \emph{strong} four coloring if there exists a four coloring of $M$ such that the regions of $M$ can be decomposed exactly into two Kempe-chains each in the form of a tree. In the light of four colorings of the maximal planar graphs given in Figure 1 and of the maps given in Figures 13,14 and 16 we assert the following conjecture.\\

\textbf{Conjecture 2}.\emph{ If the given map is hamiltonian then there exists a double-spiral ordering that results in strong four coloring of the map.}
\newpage

\begin{figure}
\centering
\includegraphics[scale=0.3, angle=-90]{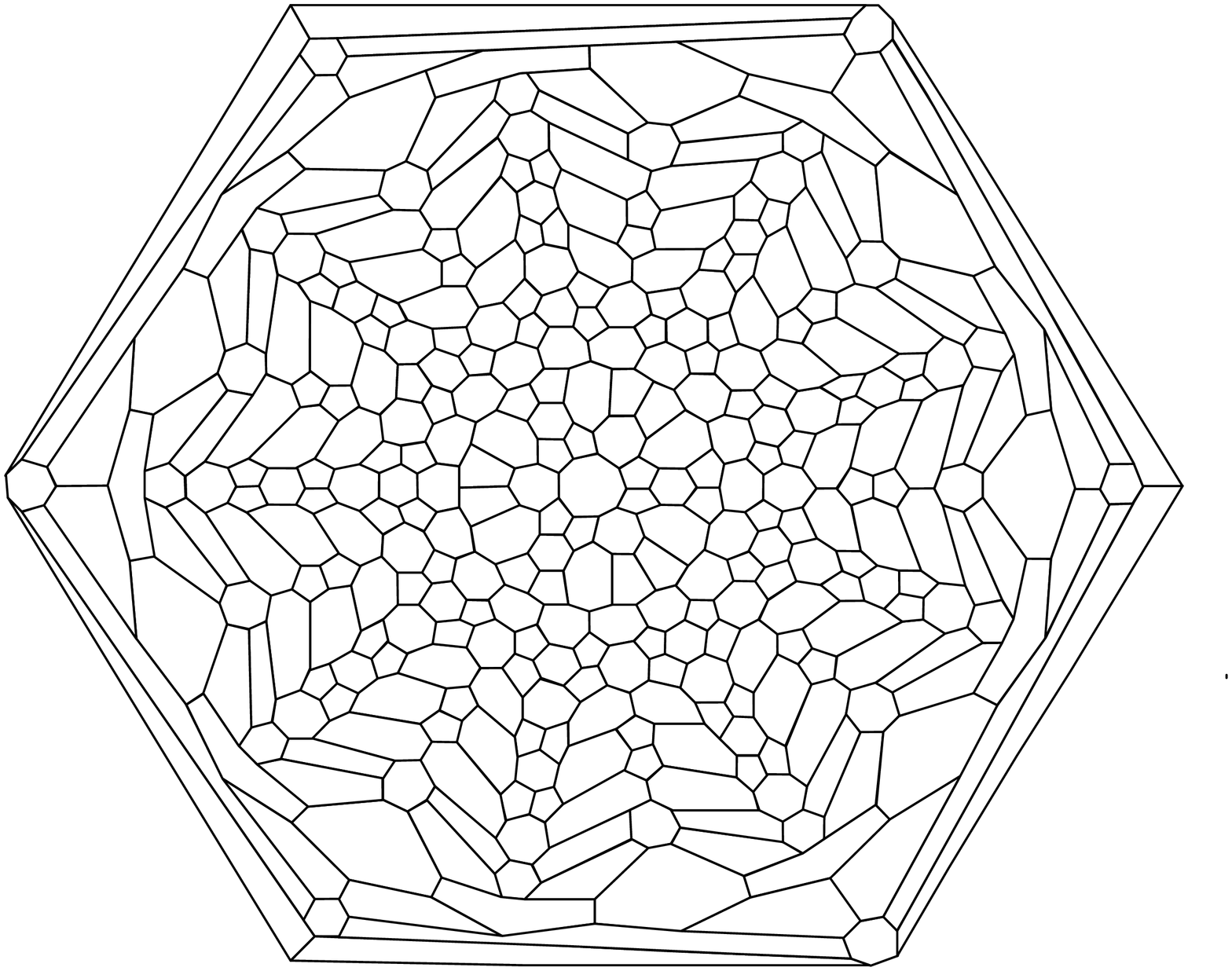}
\caption{The Haken and Appel's map. This map has been taken from Ed Pegg Jr's mathpuzzle.com/4Dec2001.htm. Haken and Appel needed a computer to $4$-color the following hardest-case map, which has been presented in a slightly different form. In this appendix we will explain step-by-step our algorithmic proof of the four color theorem on this map.}
\end{figure}

\begin{figure}
\centering
\includegraphics[scale=0.3, angle=-90]{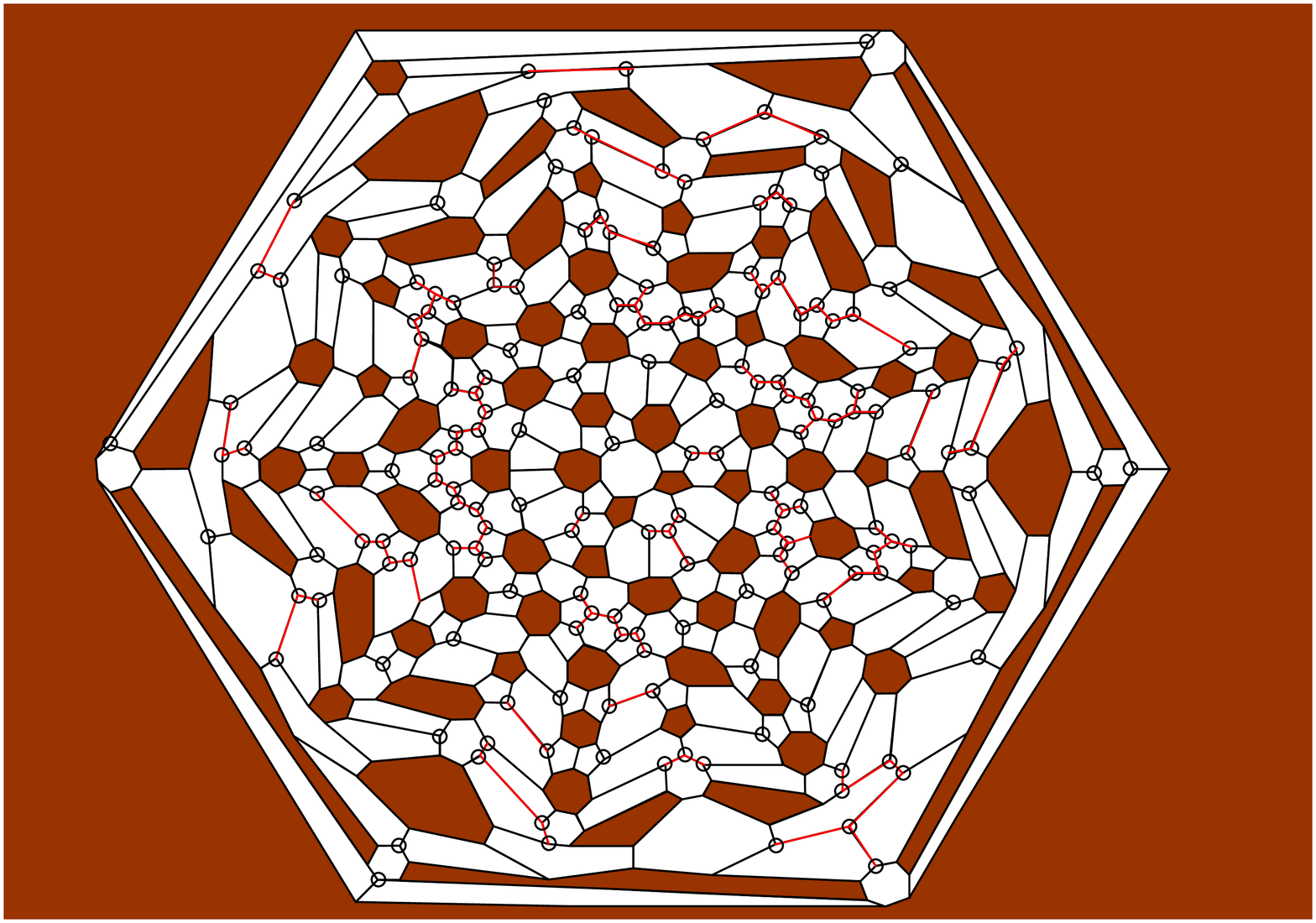}
\caption{Maximal mono-chromatic coloring of high-land (brown) regions. Note that we start coloring from the outer region and must be all adjacent to white (not colored) regions. Intersection of three adjacent regions have been shown with small circles (unwanted spots) and must be vanished as shown in Figure 10 in the maximal $2$-coloring of the map.}
\end{figure}

\begin{figure}
\centering
\includegraphics[scale=0.3, angle=-90]{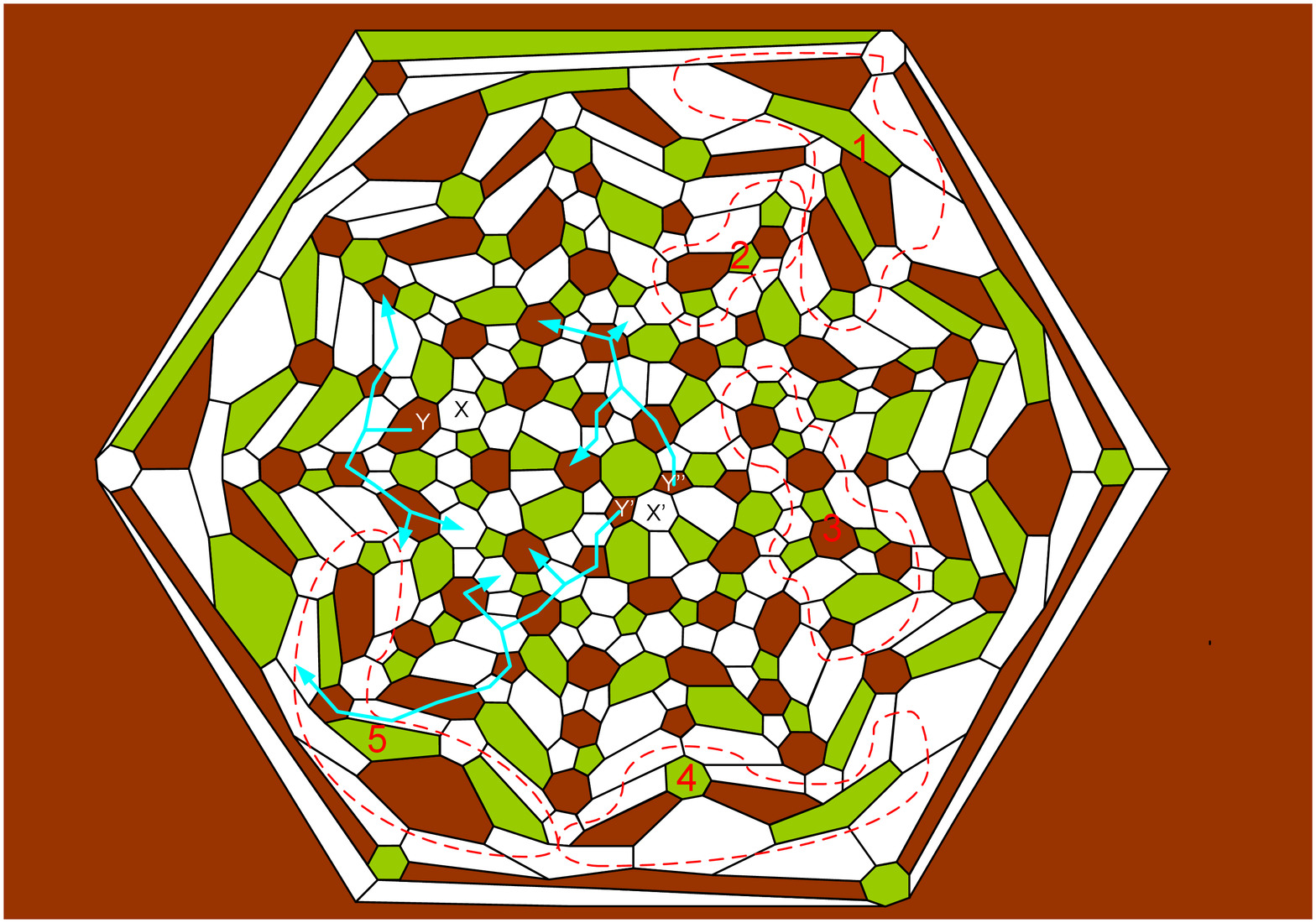}
\caption{Maximal two coloring of high-land (brown) and low-land (green) regions. Green coloring starts from the upper white region (can be started any white region adjacent to outer region). Trace of green regions form a spiraling in the clockwise direction and at each step at least one "circle" of Figure 11 is vanished by the assignment of the green color to a white region. By red-dashed curves we have shown five even white-rings (even-cycles) around the brown-green (high-lowland) islands. The rest of white regions induce an acyclic graph. The trees (blue) show $(B,W)$-Kempe chain started from region $Y$ to change the color white of the region $X$ into brown and $(B,W)$-Kempe chains started from regions $Y'$ and $Y''$ to change white $X'$ into brown. This is the only transformation that will be applied to break odd-white ring in $M(B,G)$.}
\end{figure}

\begin{figure}
\centering
\includegraphics[scale=0.3, angle=-90]{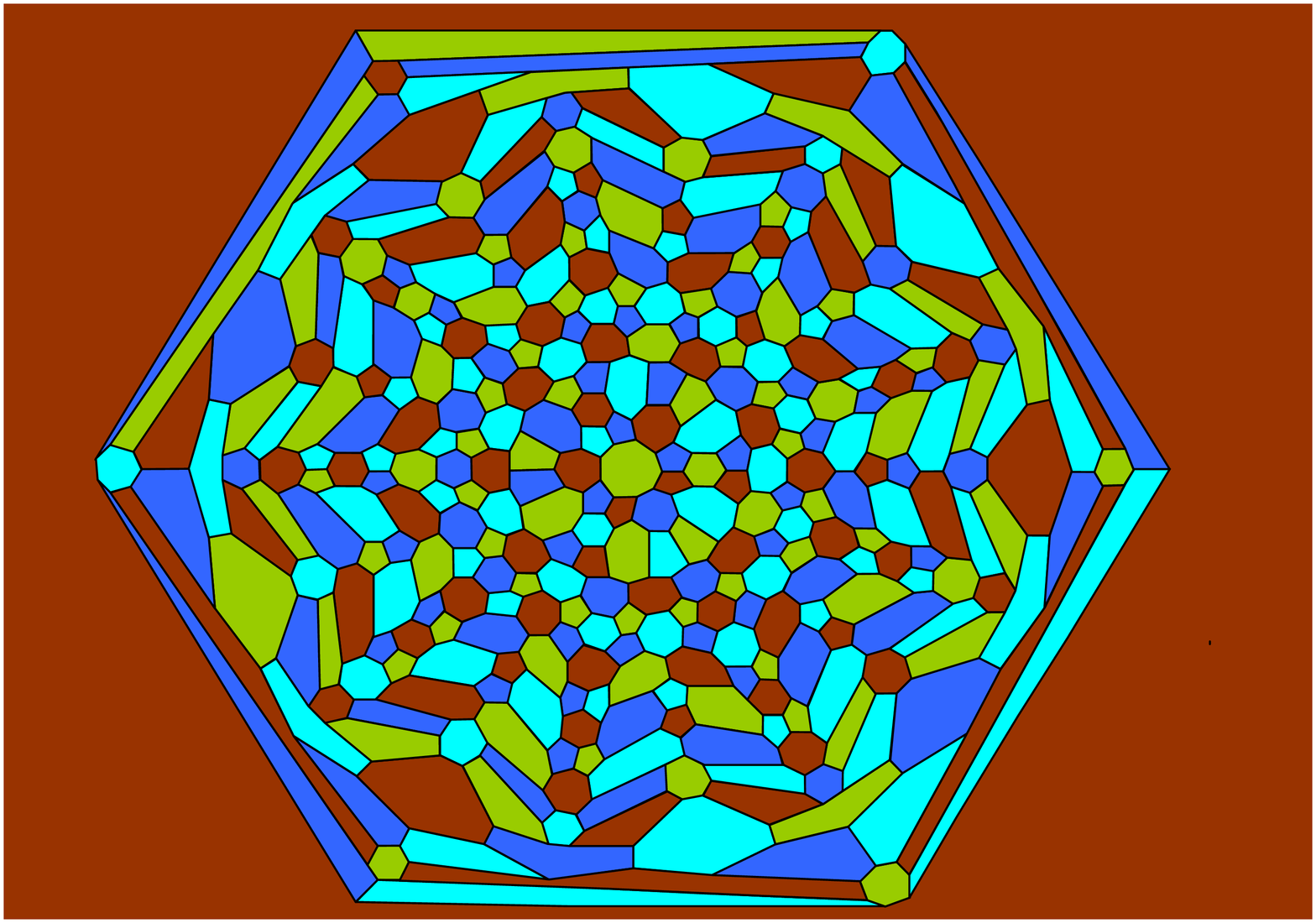}
\caption{Four coloring of Appel and Haken's map; two coloring of deep sea (dark blue) and shallow sea (light blue) regions of the two colored map of Figure 12. Here two colors is enough for the white regions since the induced dual-graph is bipartite.}
\end{figure}

\begin{figure}
\centering
\includegraphics[scale=0.55]{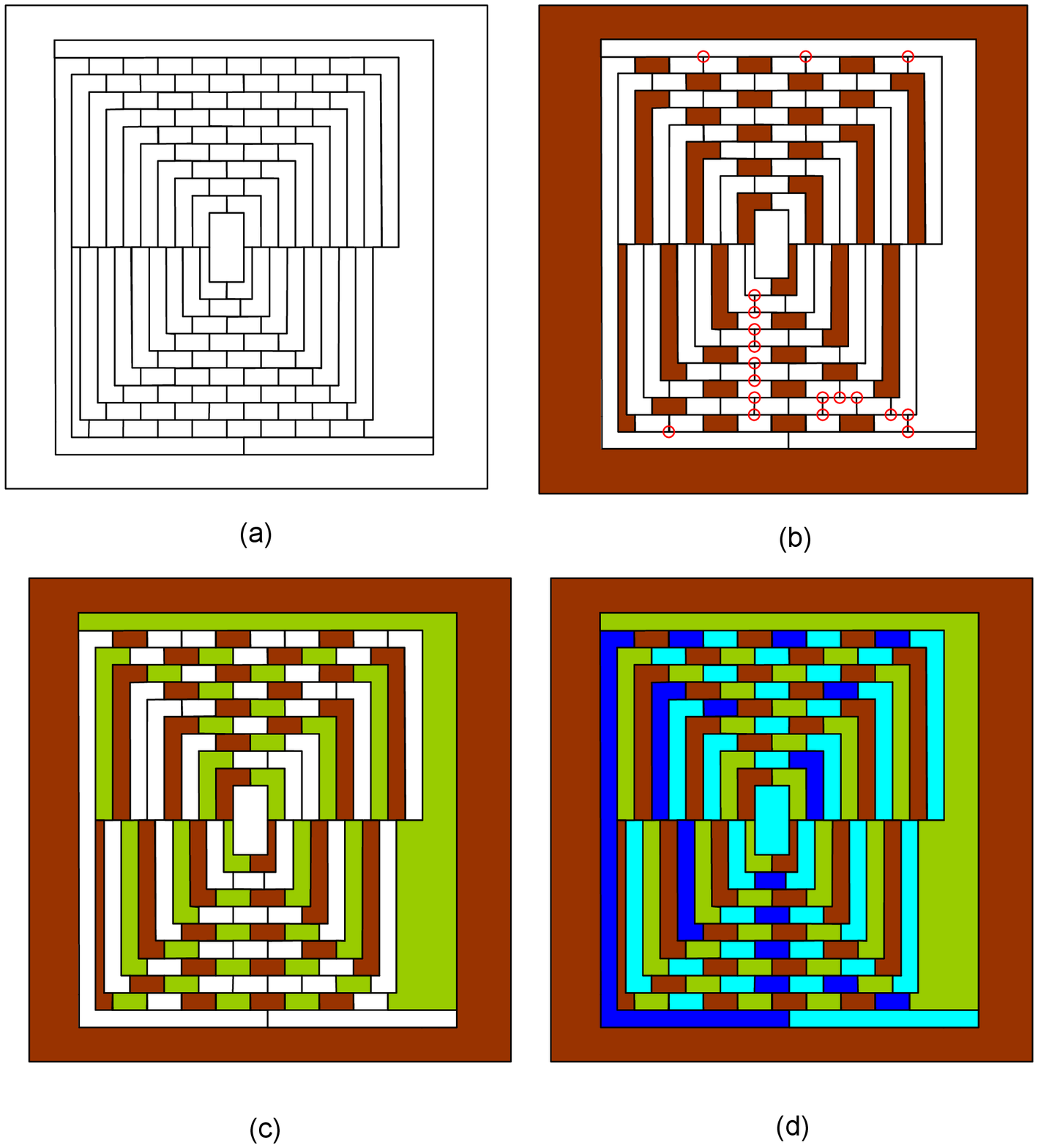}
\caption{Martin Gardner's April Fool's joke (1975). (a) The "counter-example" map, (b) Brown highland islands, (c) Brown-green high-low islands and (d) The four colored map. Note that (i) each color spiraling in the map and (ii) white regions in (c) induced disjoint union of acyclic subgraphs. Wagon has given four coloring of the April's Fool's map by using Kempe's original algorithm without facing any impasse [28],[29]. }
\end{figure}

\clearpage

\begin{figure}
\centering
\includegraphics[scale=0.4, angle=-90]{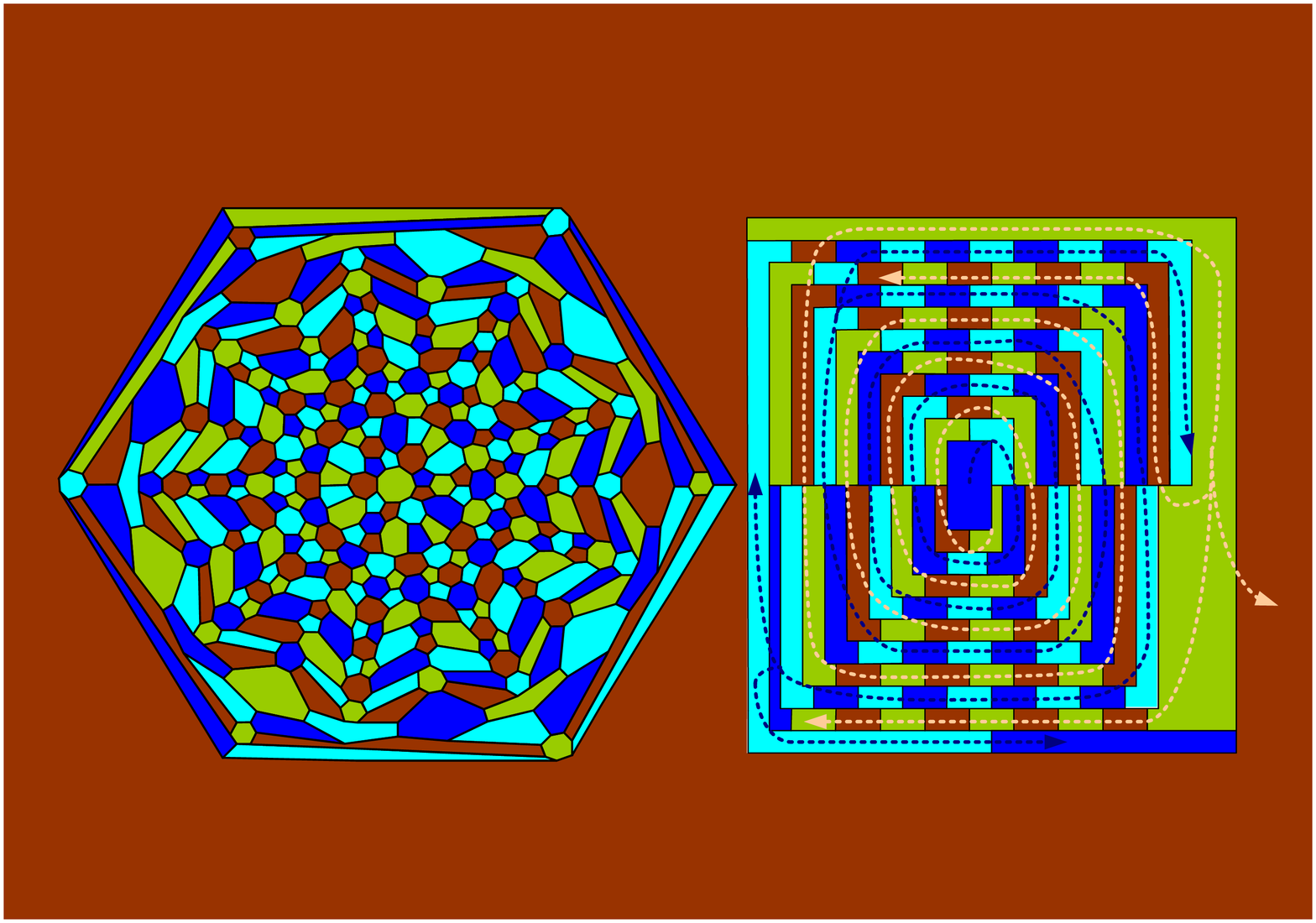}
\caption{ A normal map M has a \emph{strong} four coloring if there exists a four coloring of M such that the regions of M can be decomposed exactly into two Kempe-chains each in the form of a tree.
Stein observed that a normal map has such a coloring \emph{iff} its associated cubic graph is hamiltonian. The difficult direction of which is not difficult to see by tracing around the regions of one (hence both) of the trees. So there are many examples, though not very small ones. In the figure we have shown strong four coloring both for Appel and Haken's map and Gardner's map by the use of double spiral chain ordering and coloring. We have not proved but we suspect that if the given map is hamiltonian then double-spiral ordering results in strong four coloring of the map.}
\end{figure}

\begin{figure}
\centering
\includegraphics[scale=0.3, angle=-90]{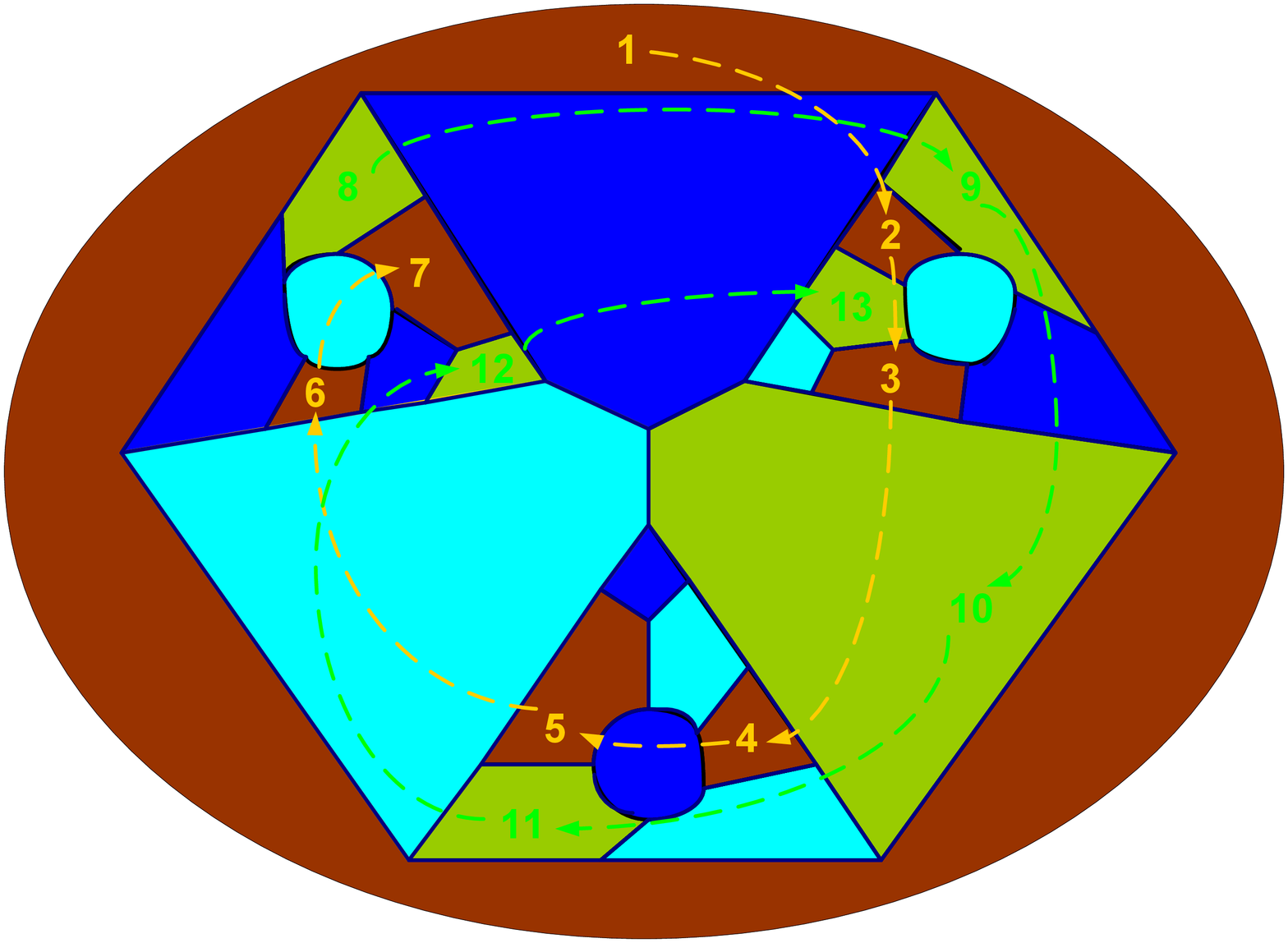}
\caption{ This is the famous Tutte's graph (1946) which disproves Tait's conjecture that every 3-regular 3-edge-connected planar graph is hamiltonian. When the Tutte's graph as seen cubic planar map (normal map) is an example of non-existence of a strong four coloring. In fact in 1971 Stein observed that a normal map has strong 4-coloring or B-set [34] iff its associated cubic graph is hamiltonian. However as shown above spiral ordering with the map coloring algorithm results a 4-colored map. }
\end{figure}

\begin{figure}
\centering
\includegraphics[scale=0.3, angle=-90]{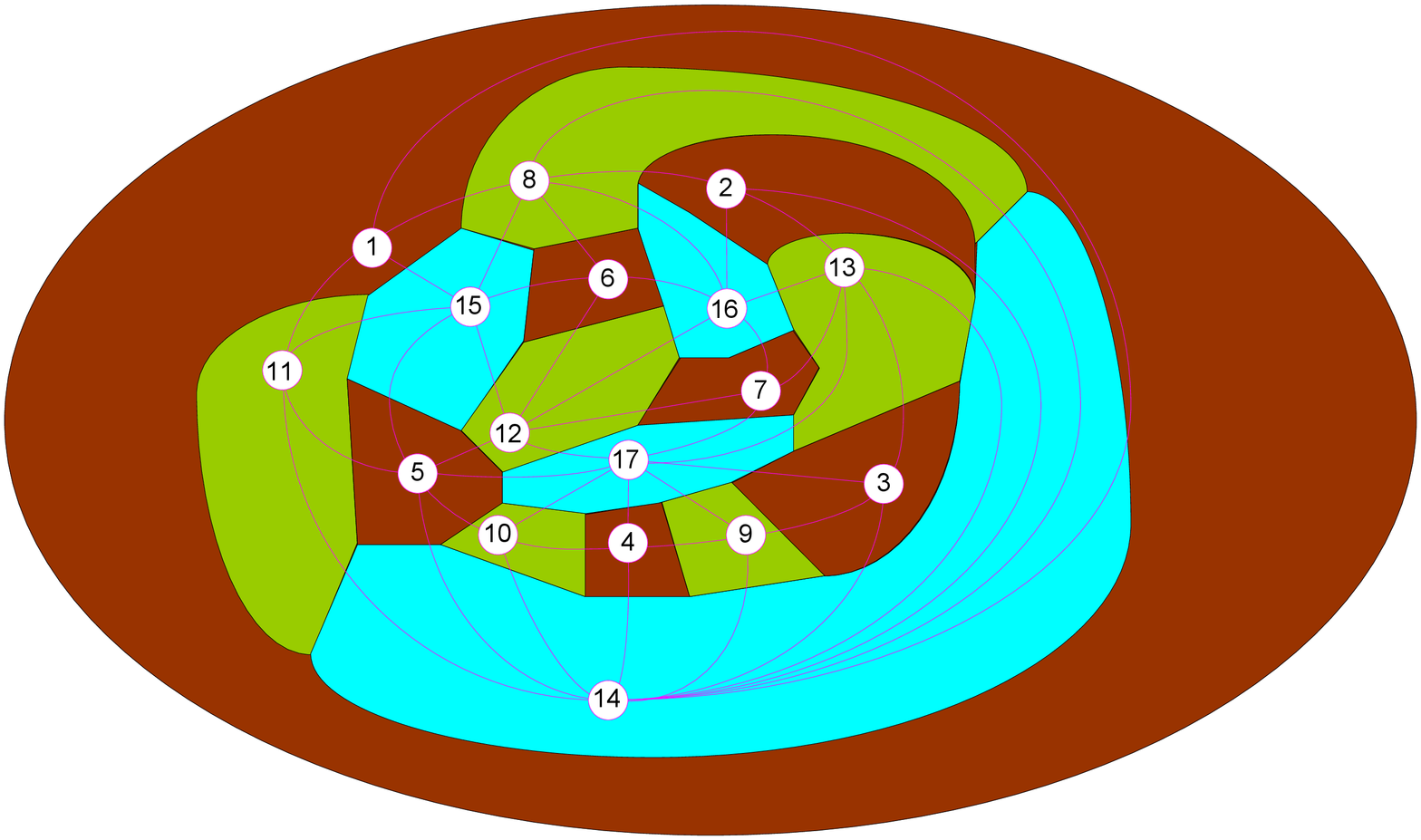}
\caption{The three colored map $M$ is an illustration of Heawood's theorem [16]:
\emph{A normal map $M$ is $3$-colorable if and only if its dual planar graph $G$ has an even triangulation.}
In the figure three coloring of $M$ is obtained by the \emph{Map Coloring Algorithm}. Numbers assigned shows the spiral order of the regions in the coloring e.g., \emph{1,2,3,4,5,6,7}  (Step 1 for high-lands), \emph{8,9,10,11,12,13} (Step 2 for low-lands)and \emph{14,15,16,17} (Step 3 for shallow-seas). Here we may call the shallow seas as  \emph{lakes} since they all surrounded by land regions. }
\end{figure}

\begin{figure}
\centering
\includegraphics[scale=0.37, angle=-90]{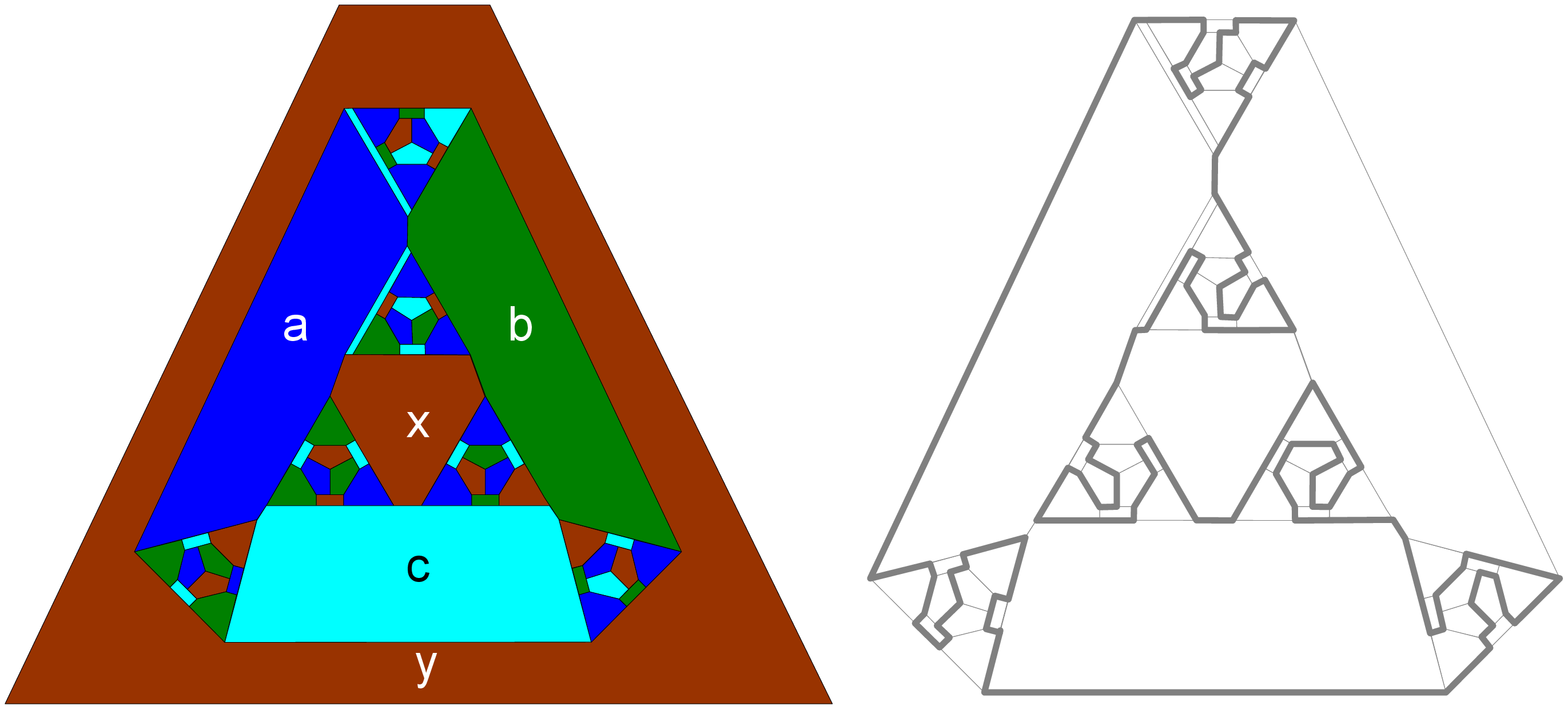}
\caption{The map shown in the figure has been provided as an counter-example map without triangle regions to the previous version of the map coloring algorithm by the referee. Since adjacency of regions labeled $x,a,b,c$ forms a complete graph $K_{4}$ on $4$ vertices, for any four coloring of the map $M$ the color of the outer-face region $y$ and central region $x$ must be same. Indeed if we do not select regions with maximum number of spots (i.e., maximum number of neighbor regions) in $M$ in Step 1 we will need the fifth color. Spiral ordering of the regions for this map results strong four coloring e.g., decomposes into two Kempe chains, hence showing that it is also Hamiltonian (see the boundaries in bold-lines in the map (right)).}
\end{figure}

\newpage
\section{Concluding remarks}
We extract the following from the first page of Appel and Haken's paper [3]:\\

\textsf{The first published attempt to prove the Four Color Theorem was made by A.B. Kempe in 1879. Kempe proved that the problem can be restricted to the consideration of "normal planar maps" in which all faces are simply connected polygons, precisely three of which meet at each vertex. For such maps he derived from Euler's formula the equation\\
\begin{center}
$4p_{2}+3p_{3}+2p_{4}+p_{5}=\sum^{k_{max}}_{k=7}(k-6)p_{k}+12$
\end{center}
where $p_{i}$ is the number of polygons with precisely $i$ neighbors and $k_{max}$ is largest value of $i$ which occurs in the map. This equation immediately implies that every maximal planar map contains polygons with fewer then six neighbors. In order to prove the Four Color Theorem by induction on the number $p$ of polygons in the map $(p=\sum p_{i})$, Kempe assumed that every normal map with $p\leq r$ is four colorable and considered a normal planar map $M_{r+1}$ with $r+1$ polygons. He distinguished the four cases that $M_{r+1}$ contained a polygon $P_{2}$ with two neighbors, or a triangle $P_{3}$ or a quadrilateral $P_{4}$, or a pentagon; at least one of these must apply by the equation.}\\

This beautiful Victorian Age deduction works for $P_{i},i=2,3,4$ and unfortunately fails for $i=5$. I think no mathematician of that period would be able to guess the possible length of a proof in future based on reducibility.\\
In this paper, by choosing direct proof, that is the opposite direction of the above, we have given an algorithmic proof for the Four Color Theorem which is based on an coloring algorithm and avoiding three-colorability in a maximal two-colorable map. The last word about the proofs given in [6],[7],[8] and including this one that uses spiral chains in the coloring algorithm. Simply enable an efficient coloring algorithm and protect us to fall in a situation similar to Kempe-tangling.

Again Appel and Haken argue strongly that [12],[13]:\\

\textsf{...it is very unlikely that one could use their proof technique without the very important aid of a computer to show that a large number of large configurations are reducible. Of course, this does not rule out the possibility of some bright \emph{young} person devising a completely new technique that would give a relatively short proof of the theorem.\\}

This paper does not prove the truth of the first sentence but it does prove that the second sentence is wrong, not only just because of the length of the proof.

\vspace{1cm}
\textbf{Acknowledgment}. The author would like to thank the referee for constructive comments; particularly for his counter-example to an earlier version of the map coloring algorithm. Thanks are also due to Tommy Jensen for discussions on the spiral ordering. Encouragements and advices given by Henry Crapo and Thomas L. Saaty are also appreciated. Finally I would like to thank to my twin brother Refik Cahit for several talks on the philosophical implications on the use of spirals in the proof of the four color theorem.

\end{document}